%% file: main.tex
\documentclass[reprint,amsmath,amssymb,aps,pra,nofootinbib ]{revtex4-2}
\usepackage{graphicx}% Include figure files
\usepackage{dcolumn}% Align table columns on decimal point
\usepackage{bm}% bold math
\usepackage{multirow}

\usepackage[linesnumbered,ruled,vlined]{algorithm2e}
\usepackage{amsmath}
\usepackage{amsfonts}

\newtheorem{theorem}{Theorem}

\begin{document}

% \preprint{APS/123-QED}

\title{A Nurse Staffing and Scheduling Problem \\ with Bounded Flexibility and Demand Uncertainty}% Force line breaks with \\
% \thanks{A footnote to the article title}%

\author{Si Zhang}
 %\email{sizhang@smu.edu.sg}
 % \altaffiliation[Also at ]{Physics Department, XYZ University.}%Lines break automatically or can be forced with \\
\author{Paul Mingzheng Tang}%
 %\email{paultang@smu.edu.sg}
\author{Hoong Chuin Lau}
 \email{Corresponding author email: hclau@smu.edu.sg}
\affiliation{%
School of Computing and Information Systems, Singapore Management University\\
Singapore
}%

% \collaboration{MUSO Collaboration}%\noaffiliation

% \author{Charlie Author}
%  \homepage{http://www.Second.institution.edu/~Charlie.Author}
% \affiliation{
%  Second institution and/or address\\
%  This line break forced% with \\
% }%
% \affiliation{
%  Third institution, the second for Charlie Author
% }%
% \author{Delta Author}
% \affiliation{%
%  Authors' institution and/or address\\
%  This line break forced with \textbackslash\textbackslash
% }%

% \collaboration{CLEO Collaboration}%\noaffiliation

\date{May 28, 2025}% It is always \today, today,
             %  but any date may be explicitly specified

\begin{abstract}
Nurse staffing and scheduling are persistent challenges in healthcare due to demand fluctuations and individual nurse preferences. This study introduces the concept of bounded flexibility, balancing nurse satisfaction with strict rostering rules, particularly a real-world time regularity policy from a major hospital in Singapore. We model the problem as a multi-stage stochastic program to address evolving demand, optimizing both aggregate staffing and detailed scheduling decisions. A reformulation into a two-stage structure using block-separable recourse reduces computational burden without loss of accuracy. To solve the problem efficiently, we develop a Generative AI-guided algorithm. Numerical experiments with real hospital data show substantial cost savings and improved nurse flexibility with minimal compromise to schedule regularity. Numerical experiments based on real-world nurse profiles, nurse preferences, and patient demand data are conducted to evaluate the performance of the proposed methods. Our results demonstrate that the stochastic model achieves significant cost savings compared to the deterministic model. Notably, a slight reduction in the regularity level can remarkably enhance nurse flexibility.

% \begin{description}
% \item[Usage]
% Secondary publications and information retrieval purposes.
% \item[Structure]
% You may use the \texttt{description} environment to structure your abstract;
% use the optional argument of the \verb+\item+ command to give the category of each item. 
% \end{description}
\end{abstract}

%\keywords{Suggested keywords}%Use showkeys class option if keyword
                              %display desired
\maketitle

%\tableofcontents

\input{01_intro}
\input{02_lit_rev}
\input{03_model_form}
\input{04_sol_approach}
\input{05_experiments}

\input{06_conclusion}

\section{Acknowledgments}
This research is supported by the National Research Foundation, Singapore under its AI Singapore Programme (AISG Award No: AISG2-AISG2-100E-2023-118). The authors thank the Nursing Department of Tan Tock Seng Hospital for providing the problem statement on which our model is built. 

% The \nocite command causes all entries in a bibliography to be printed out
% whether or not they are actually referenced in the text. This is appropriate
% for the sample file to show the different styles of references, but authors
% most likely will not want to use it.
\nocite{*}

\bibliography{main}% Produces the bibliography via BibTeX.

\input{07_appendix}

\end{document}

%% file: 01_intro.tex
\section{Introduction}
\label{Section_Introduction}

Nursing resources constitute a significant proportion of caregivers in healthcare systems worldwide. However, the nursing profession continues to face shortages due to ever-increasing healthcare demands \citep{halter2017determinants, Ringo2021}. These shortages are attributed to several factors, including an aging population, lack of flexibility, and irregular work patterns \citep{breugem2022equality}. As a result, healthcare managers must make strategic nurse planning and operational decisions to ensure highly efficient labor utilization. To address this challenge, this paper investigates a nurse staffing and scheduling problem based on a real-world application scenario in a hospital in Singapore.

The hospital we studied faces several challenges in its nurse management system. First, patient demand is uncertain and often unknown in advance. For example, in practice, accurate demand information is typically available only one week earlier \citep{ceschia2019second}. Ignoring these uncertainties can lead to suboptimal nursing schedules. Second, feasible schedules are constrained by diverse rostering rules, including supply-demand balance, workload distribution, pattern requirements, and fairness considerations, which further increase the complexity of the problem. Third, individual nurses increasingly seek work-time flexibility, such as accommodating work-time preferences and specific requests. Consequently, to enhance staff morale and improve retention, nurse scheduling systems must account for individual flexibility as much as possible. However, this often conflicts with the hospital's rostering rules from the management perspective. One of the most contentious rules in this study is based on a real-world hospital policy that enforces a specific level of time regularity in nurses' schedules, which should not be violated. To address this, we propose a new concept, termed bounded flexibility, to capture the optimal balance between satisfying nurses' preferences and adhering to managerial rostering rules.

In this paper, we address a nurse staffing and scheduling problem that incorporates our proposed concept of bounded flexibility and accounts for demand uncertainty. Within a given planning horizon, staffing decisions aim to optimize the number of nurses available for scheduling, while scheduling decisions focus on determining detailed shift arrangements for the nursing staff. To ensure adequate coverage of uncertain demand, we develop an integrated staffing and scheduling model that generates an initial and staffing plan and schedule (here and now solution), while providing for subsequent adjustments (wait and see solutions) in response to actual demand as it becomes known stage-by-stage. The proposed formulation distinguishes clearly between hard rules, which must be strictly adhered to, and flexible soft rules, which enable bounded flexibility to accommodate nurse preferences.

To summarize, we make the following contributions in this paper. 
First, a deterministic model is developed to optimize the integrated nurse staffing and scheduling decisions under the perfect information of patient demand. The aim is to minimize the sum of staffing cost, demand mismatch penalty, requests rejection penalty, and pattern violation penalty. 
Unlike the traditional nurse scheduling problem \citep{miller1976nurse,warner1976scheduling,aickelin2004indirect}, the integrated staffing and scheduling model enhances overall caregiving efficiency by right-sizing staff. In this model, nurse staffing plans and scheduling decisions are interdependent, and optimizing them independently often results in suboptimal outcomes, ultimately compromising the quality of care.
Second, given all compulsory rostering rules from the hospital level, bounded flexibility is introduced to help accommodate nurse-specific needs as much as possible. Imposing flexibility into nurse schedules that satisfy the preference of each nurse, while considering hospital management policies, nurse workload, and government regulations, is the key to ensuring nurse welfare. 
Third, a multi-stage stochastic programming model is developed to incorporate the evolving demand over time by dividing the planning horizon into several stages, following the approach taken by \citep{Kim2015}. Our proposed model determines an initial schedule, apart from which two levels of decisions (i.e., aggregate-level and detailed-level) for each stage are optimized subsequently for each demand scenario. Specifically, the aggregate-level decision with nurse staffing is generated before the actual demand realization of this stage. After demand is revealed, the detailed-level decision is made to determine an operational nurse schedule plan allowing short-term adjustments to the initial schedule. Note that these adjustments are done without having to re-solve the model, since all solutions have been precomputed by the model.  Broadening the problem scope to multi-stage substantially increases the solving difficulty. Therefore, the special structural property for the model is identified and the multi-stage program with \textit{block-separable recourse} described by \cite{birge2011introduction} can be equivalently reduced to a two-stage stochastic program. This property allows us to significantly shorten the solution time. 
Fourth, a state-of-the-art generative approach based on Generative Flow Networks (GFN, \citet{Bengio2021a,Bengio2021b}) is developed to obtain near-optimal solutions within reasonable computational times. This approach concurrently generates a diverse set of candidate nurse schedules by training a policy to sample schedules in proportion to their rewards. Once the actual demand is realized, the best-performing schedule is selected from this set of candidates. To evaluate the model's performance, a real-world case study was conducted at a surgical ward at Tan Tock Seng Hospital in Singapore \citep{TTSH}, utilizing nurse profile data, nurse preference data, and patient demand data. \footnote{Due to data sensitivity, the results reported in this paper are based on anonymized data, and no private information is disclosed.}
Our findings offer several managerial insights. First, the model demonstrates significant cost savings by accounting for demand uncertainty and its evolving nature. Additionally, even a small reduction in the regularity level can substantially increase nurse flexibility.

The remainder of the paper is organized as follows: Section ~\ref{Section-Literature} reviews the relevant literature. Section ~\ref{Section-model formulation} presents the detailed mathematical model. Section ~\ref{Section-solution method} describes the proposed solution approach. Section ~\ref{Section-empirical study} reports the results of a case study conducted in a hospital in Singapore. Finally, Section ~\ref{Section-conclusions} summarizes the conclusions and highlights future research directions.

%% file: 02_lit_rev.tex
\section{Literature review}
\label{Section-Literature}
In this section, two streams of studies relevant to this paper are reviewed. Section ~\ref{Review-NSP} includes the studies on nurse scheduling problems. Section ~\ref{Review-GFN} gives a summary of studies on generative AI models for solving scheduling problems. In addition, a summary of the comparison of the relevant papers along different dimensions (demand settings, fairness, flexibility, etc) is provided in Appendix \ref{appen-studies}.

\subsection{Nurse scheduling problem}
\label{Review-NSP}
Nurse scheduling problems are complex due to inherent uncertainties of patient demand. An extensive review on healthcare planning and scheduling over the past 30 years is in \cite{youn2022planning}. \cite{Ricardo2023} proposes a preventive and reactive approach to study unexpected nurse absence. A multi-objective nurse scheduling model is presented that determines the backup team and base schedule for each day and two reactive rescheduling heuristics are further given. \cite{Schoenfelder2020} incorporates two quick-response approaches with patient demand fluctuations, i.e., (i) allowing the transfer of float nurses across units and (ii) allowing the transfer of patients between units and off-unit admissions. \cite{Legrain2020} studies a personalized nurse scheduling problem under uncertainty.  Before the realization of future inputs, for the current stage, a dynamic approach is presented to make an irrevocable schedule, with weekly patient demand and nurses' preferences as inputs.
\cite{Kim2015} investigates the integrated staffing and nurse scheduling problem with stochastic patient demand. They establish a two-stage stochastic programming model and develop a modified multi-cut approach based on an integer L-shaped algorithm for solving the model.  \cite{Maenhout2013b} study the nurse re-rostering problem with unexpected events from nurses. An optimization/heuristic tool is used that makes adjustments to the schedule of heterogeneous nurses, together with a branch-and-price algorithm or a meta-heuristic optimization. 

Some studies are designed with preference considerations, like preferred/disliked shift types, working days, or specific patterns \citep{Bard2005,Bard2007}. Considering that nurses may treat various shift types differently, \cite{Ricardo2023} propose a method to model nurses' preferences of shifts by scaling them from 1 to 5 as marks. A bi-objective model is proposed maximizing the total satisfaction of nurses. In order to minimize the uncovered patients and violated nurses' preferences, \cite{Guo2022} proposes a mixed-integer programming model for nurse scheduling problems,  the sum of weighted uncovered demand and nurse preference violations over one month. Several undesired patterns are considered, including assignments of different shift types on two consecutive days and nonworking, working, nonworking successively on three consecutive days.

A handful of studies consider dynamic effects in the nurse scheduling problem, which means that the decisions are taken multiple times or in multiple stages, with the revelation of information. Considering the dynamic feature, \cite{Kheiri2021} presents a multi-stage nurse scheduling model. A sequence-based selection hyper-heuristic utilizing a statistical Markov model is introduced.  \cite{Schoenfelder2020} formulates a multi-stage stochastic program, which determines the initial nurse-to-shift assignment, and the following nurse-to-unit assignment. The objective is to minimize the costs of regular and overtime shift assignments, the expected understaffing penalties, the patient turnaways penalties and the patient transfer penalties. Considering the weekly given manner of patient demand and preferences of nurses, \cite{Legrain2020} proposes an online stochastic algorithm with sample average approximation. First, the method is able to generate a small set of candidate schedules. Then, with various test scenarios, they rank all candidates with an evaluation process, after which keep the best-performing one.

\subsection{Generative AI models}
\label{Review-GFN}
Solution methodologies for nurse scheduling problems can be mainly categorized as meta-heuristics, exact algorithms, and hybrid approaches.  In recent years, learning-based methods have been introduced to enhance the efficiency of optimization-based methods in staff scheduling problems. In order to automate the design of heuristics and reduce complexity, \cite{Chen2022} proposes an approach assisted with a neural network. The selection of heuristic is under the guidance of a deep neural network model, and they use a recurrent neural network model to avoid local optimal solutions.  For a complex staff re-rostering problem, \cite{Oberweger2022} proposes an enhanced large neighborhood search method. A novel destroy operator is designed, which uses a conditional generative model to identify variables to be selected and then integrates with a designed sampling strategy to make actual selections.   \cite{Park2021} presents a solution method for solving a job-shop scheduling problem with a graph neural network and reinforcement learning. Considering the feature of the problem, job-shop scheduling is formulated as a sequential decision-making process, with a graph to describe its states. 

Among learning-based methodologies, generative models are designed to learn the underlying distributions of the data, versus discriminative models which learn the decision boundaries of data. Generative flow networks are models that can be used to generate a solution or object by learning a stochastic policy function. The probability of generating a particular solution/object is proportional to the reward function. GFN has been extensively used in different fields, including scheduling \citep{zhang2023robust}, transportation safety \citep{Zhang2023}, and marketing \citep{Liu2023}. Relying on the efficiency of modeling a scheduling process with a directed acyclic graph, \cite{zhang2023robust} proposes a generative scheduling approach, that aims at training a generative model that generates a set of candidate schedules that have a low makespan and are diverse.  GFN is introduced as the method for learning the stochastic policy that generates the scheduling graph piece-by-piece and is proportional to a given reward.

%% file: 03_model_form.tex
\section{Model formulation}
\label{Section-model formulation}
This section presents two models. We first present a deterministic model with a full set of rostering rules, assuming that patient demand is known apriori, which is given in Section \ref{Section-determi model}. Then, in Section \ref{Section-stoch model}, a multi-stage stochastic programming model is formulated given uncertain demand realization from stage to stage.

\subsection{Deterministic model}
\label{Section-determi model}

\begin{table*}
\caption{\textbf{Rostering rules}}
\centering
        {\def\arraystretch{1.2}  %change the number for increasing or decreasing the spacing.
            \scriptsize %uncomment this for changing the font size
\begin{tabular}{llll}
\hline\hline
\textbf{Categories}                                     & \multicolumn{2}{l}{\textbf{Rules}} & \textbf{Descriptions} \\ \hline
\multirow{4}{*}{\textbf{Supply and demand rules}} & \multirow{2}{*}{Supply}   & {H1. Nurse time preferences}  &  Nurses can only be assigned shifts within their preferred times. \\
                                               &                           & {S1. Nurse requests}  & The specific leave or shift requests of nurses should be \\
                                               & \multirow{2}{*}{}  & & considered.            \\ \cline{2-4}
                                               & \multirow{1}{*}{Demand}   & {S2. Patient demand}  & For each day and each time slot, a certain number of nurses are             \\
                                               & \multirow{1}{*}{}   &  & required  to cover the patient demand.             \\
                                               \hline
\multirow{6}{*}{\textbf{Workload rules}}          & \multicolumn{2}{l}{{H2. Work hours}}          &For each nurse, the overall work hours in the planning horizon              \\
&\multicolumn{2}{l}{}   &must be within a range according to the contract.              \\
                                               & \multicolumn{2}{l}{{H3. Single assignment}}          &   Each nurse can only be assigned one shift per day.           \\
                                               & \multicolumn{2}{l}{{H4. Work policy}}          & Each nurse has a work policy restricting the schedule regularity.        \\
                                               % &\multicolumn{2}{l}{}   & schedule.          \\
                                               % & \multicolumn{2}{l}{{H5. Rest days per week}}          &   Each nurse must have certain rest days per week.           \\
                                               % % & \multicolumn{2}{l}{{H6. Distance of N cluster}}          & For each nurse, the assignment of two adjacent N clusters must be \\
                                               % % &\multicolumn{2}{l}{}   & separated by enough rest days. N cluster refers to consecutive night shifts.        \\
                                               & \multicolumn{2}{l}{{S3. Maximum W weekends}}          & Each nurse is desired to be assigned a minimum of certain              \\
                                               &\multicolumn{2}{l}{}   & weekend days as rests day.       \\\hline
\multirow{2}{*}{\textbf{Pattern rules}}           & \multicolumn{2}{l}{{H5. Ban patterns}}          & Some patterns must be avoided on consecutive days, e.g., nurses        \\
& \multicolumn{2}{l}{}          & can not have an AM shift right after a day with N shift.      \\
                                               & \multicolumn{2}{l}{{S4. Undesired pattern}}         & Some patterns are undesired on consecutive days, e.g., nurses are             \\
                                               & \multicolumn{2}{l}{}          & unwilling to have a NW before and after a W.       \\
                                               \hline
\multirow{1}{*}{\textbf{Fairness rules}}          & \multicolumn{2}{l}{{H6. Rest days per stage} }         & Each nurse must have a certain equal number of rest days per stage. \\
 % & \multicolumn{2}{l}{}         & horizon should be within a certain range. \\
                                               % & \multicolumn{2}{l}{H13}         &              \\
\hline
\hline
\end{tabular}
}
\label{table-constraints}
\end{table*}

\begin{figure}
        \centering
    \caption{An illustrative example of work policy}
 \includegraphics[trim = 110 75 150 0, width=0.9 \linewidth]{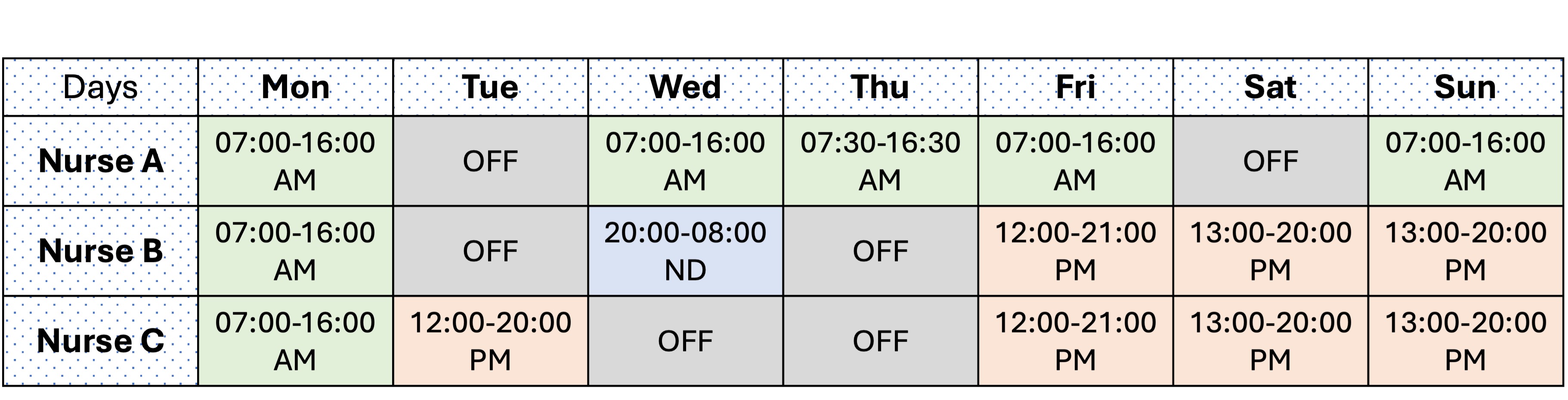}       
 \label{Fig-policy}
\end{figure}

Consider the staffing and scheduling of a given group of nurses $i\in\mathcal{I}$ within a ward over a planning horizon $\mathcal{D}$, which is discretized into several equal and short intervals $d \in \mathcal{D}$, i.e., one day. Usually, the planning horizon contains several weeks. On the hospital level, the schedule is restricted by a set of hard and soft rostering rules (Table~\ref{table-constraints}) categorized into four groups. The hard rules (H) determine the feasible region for the solution and the soft rules (S) are associated with violation penalties. Each day $d$ is further divided into three \textit{time slots} $j\in\mathcal{J}=\{AM, PM, N\}$ \footnote{For conciseness, refer to the Appendix \ref{app-abbrev} for the main abbreviations.}.  A \textit{shift} $s$ is a period with a beginning and an ending time, i.e., 8:00-16:00. Each shift belongs to one and only one of the time slot $j$, according to the time period. An \textit{assignment} is a (nurse $i$, day $d$, shift $s$) tuple. Patient demand $Q_{j}^{d}$ specifies the total number of nurses needed for a day $d$ on a time slot $j$. For each nurse, the overall roster is made of up several assignments.  Furthermore, each nurse has a predetermined \textit{work policy} $p\in\mathcal{P}=\{p_{1}, p_{2}, p_{3}\}$ according to the contract that stipulates the time regularity of their schedule within the planning horizon. To be specific, a work policy for a nurse defines the number of different time slots that can be assigned in the overall schedule. For example, in Figure \ref{Fig-policy}, we give an illustrative example of three different work policies, assuming that the planning horizon is one week. Nurse A is a $p_1$ nurse with all shifts in the AM slot. Nurse B's work policy is $p_{3}$ so this nurse has been assigned shifts that belong to all of the three time slots. Nurse C has the work policy $p_{2}$ and only works on AM and PM slots. Obviously, Nurse A has the most regular work schedule, while Nurse B tends to have the most "challenging" schedule with alternating day and night shifts, which are undesirable. Note that work policy is a hospital-level hard rule and cannot be violated. While each individual nurse has an emerging need for flexibility, in two aspects. First, time preferences. Each nurse $i$ has a subset of preferred shifts $\mathcal{S}_{i}\subset \mathcal{S}$. Second, requests. Each nurse can request some specific assignments or leave days, e.g., a night shift on Friday, taking leave on Saturday, etc (for simplicity we will term all these as requests). However, nurse flexibility often conflicts with specific rostering rules at the hospital level. In this paper, we consider bounded flexibility by satisfying nurse flexibility as much as possible while ensuring the management rules. The objective is to minimize the total cost across the planning horizon, including the staffing cost, the penalty for demand mismatch, the penalty for nurse request rejection, and the penalty for soft rule violations. The decision variables are nurse staffing ($\beta_{i}$, if nurse $i$ is scheduled) and scheduling ($x_{i}^{d,s}$, if nurse $i$ is assigned to shift $s$ on day $d$). 

Notations for the deterministic model are given in Appendix \ref{app-notations-deterministic}. The full deterministic model is formulated as shown in Model 1 in Appendix \ref{app-model1} (Equations \ref{model1-obj} to \ref{model1-variable domain_p3}).

Objective function in Equation ~\ref{model1-obj} sums the total staffing cost, the penalties of supply-demand mismatch (S2), unsatisfied nurse requests(S1) and other undesired violations. $v_{i}^{m}$ is a binary variable accounting for the total number of violations for nurse $i$. Several violations are considered, including maximum W weekend days (S3), and undesired pattern violations (S4), i.e., NW-AM, N-N on weekends and NW-AM-NW.

Constraint ~\ref{model1-num of nurses2} assures the number of scheduled nurses can not exceed the total number of available nurses in the ward. Constraint ~\ref{model1-h1} enforces that a nurse per day can only be assigned at most one shift. Note that here we enforce the nurse flexibility of work time preferences by $s\in\mathcal{S}_{i}$. In real-world applications, before the scheduling period, the nurses are allowed to submit their preferred work shifts for the entire planning horizon.    Constraint ~\ref{model1-h2} sets the upper and lower limit on the total number of work hours across the planning horizon for each nurse. Constraints ~\ref{model1-h5}-\ref{model1-h6} jointly define the work policy. Specifically, Constraint ~\ref{model1-h5} ensures that the number of different time slots scheduled for nurse $i$ is less than $W_{p}$ according to the work policy. Constraint ~\ref{model1-h6} guarantees that only if time slot $j$ is assigned to nurse $i$, then the shifts that belong to that time slot can be assigned to the nurse. To ensure enough rest hours, Constraint ~\ref{model1-h7} excludes a working pattern for nurses, that is N-NW-AM on three consecutive days. Constraint ~\ref{model1-h7b} ensures that a nurse can not be assigned to shift in the AM or PM slot right after N. Constraint ~\ref{model1-h8} restricts that there should be at least $E_{i}$ NW days per week for nurse $i$. Constraint ~\ref{model1-h16} limits the number of consecutive working days to be no more than $F_{i}$  for nurse $i$.  Constraint ~\ref{model1-h13} expresses the under- and over-supplied demand.

Also, a set of soft rule violations is specified by Constraints ~\ref{model1-minNWweekends}-\ref{model1-violations}, where intermediate variables are used to count the number of violations accordingly. In Constraint ~\ref{model1-minNWweekends}, integer variable $p_{1,i}$ counts the number of working weekend days more than the upper bound $M_{i}$, where $M_{i}$ represents the maximum number of working weekend days in the planning horizon for nurse $i$. Then some undesired pattern violations are given. In Constraint ~\ref{model1-h14}, when there is an undesired pattern NW-AM from day $d$ of nurse $i$, the binary variable $p_{2,i}^{d}=1$, indicating one violation.  Constraint ~\ref{model1-h17} calculates the violations of two consecutive N slot shifts on weekends, where $\mathcal{D}^{Sat}$ is the subset of Saturday in the planning horizon. So, if nurse $i$ has been assigned two N slot shifts on weekends from day $d$, binary variable $p_{3,i}^{d}$ is 1. Similarly, a binary variable $p_{4,i}^{d}$ is used to represent the "0-1-0" pattern in Constraint~\ref{model1-h18}.  That is to say, if a nurse $i$ is NW on day $d$ and $d+2$, and W on day $d+1$,  then $p_{4,i}^{d}$ equals to 1. Constraint ~\ref{model1-violations} calculates the total number of violations in the planning horizon for each nurse.  Finally, constraints ~\ref{model1-variable domain_xy} - \ref{model1-variable domain_p3} define the domains of the variables, where $\mathbb{Z}^{+}$ represents the set of non-negative integer variables.

Note that the above model is nonlinear due to the term of the requests in the objective function ~\ref{model1-obj}. By introducing an auxiliary binary variable $Y_{i}^{d,s}$, the model can be easily linearized. The detailed linear model (Model 2) is given in Appendix \ref{app-linear-deterministic}. Henceforth, to simplify model expressions, we use bold vectors to denote decision variables, e.g., $\mathbf{x} = \left[ {x}_{i}^{d,s}, i \in \mathcal{I}, d \in \mathcal{D}, s \in \mathcal{S} \right]$. Then, let $\boldsymbol{\chi}=\left\{\boldsymbol{\beta},\mathbf{x,y,\bar{q},\tilde{q}}\right\}$ represent the set of the vectors for all main variables, and omit the induced and auxiliary variables.

\subsection{Stochastic model}
\label{Section-stoch model}
In response to patient demand uncertainty, a multi-stage stochastic programming model is established in this section. At the beginning of the planning horizon, the model determines an initial staffing and scheduling plan. That is, it optimizes the number of nurses available to be scheduled in the planning horizon. The initial scheduling decision generates a detailed schedule plan for each nurse, before knowing the precise patient demand. The planning horizon is divided into multiple stages $h=1,2,...|\mathcal{H}|$ and $h\in\mathcal{H}$. Then, with the realization of uncertain information stage by stage, short-term adjustments of schedules are required to ensure the coverage of stochastic demand. The objective is to minimize the summation of initial staffing and scheduling costs, and the expected recourse of all stages. To facilitate modeling, we make the following assumptions.

\subsubsection{Decision Groups}
\label{Assumption-groups of decisions}
\begin{figure*}
        \centering
    \caption{An illustrative example of multi-stage decision process}
 \includegraphics[trim = 0 0 0 -10, scale=0.49,clip]{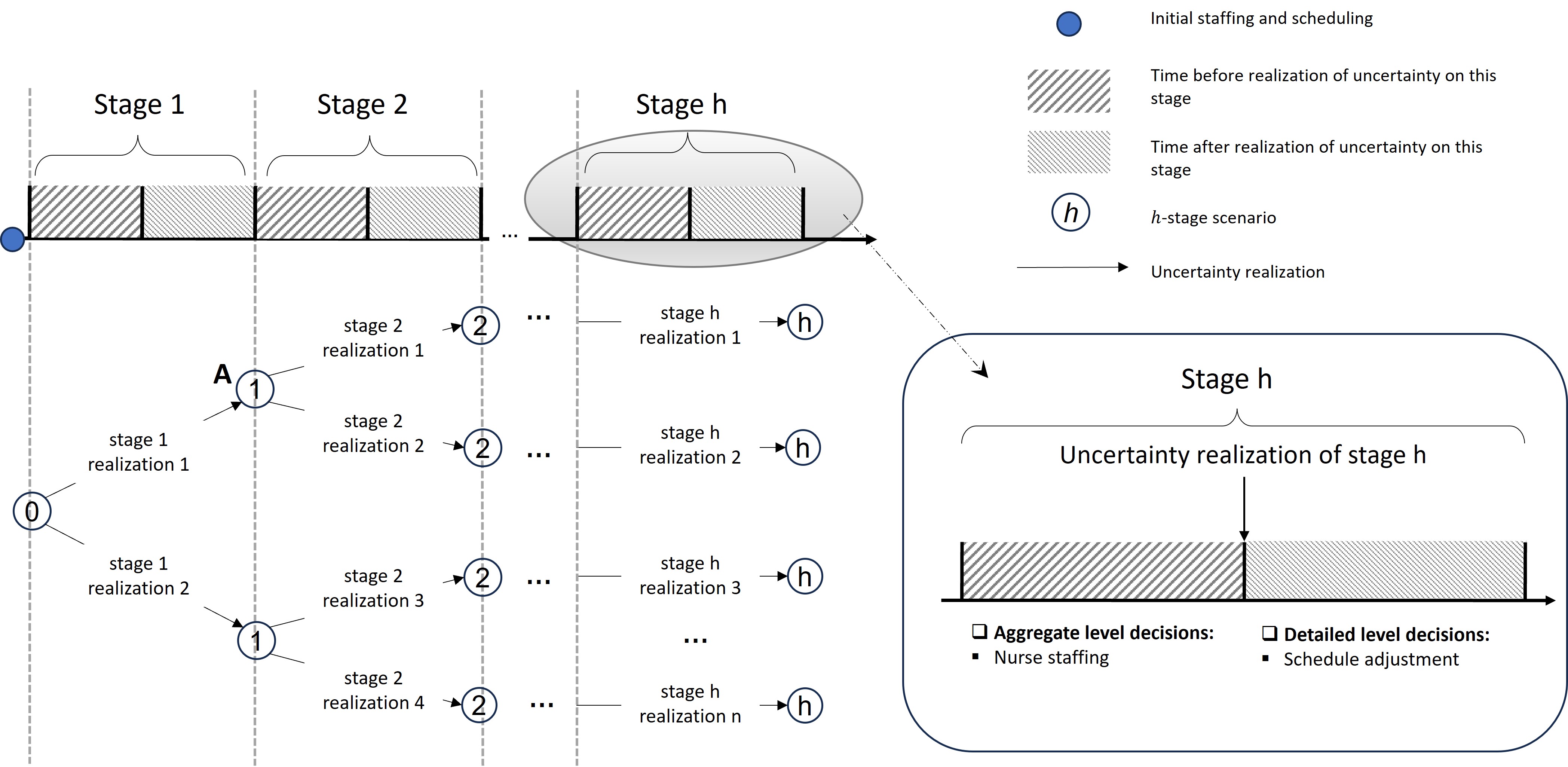}\label{Fig-multi-stage procedure}
\end{figure*}

Within each stage $h$, there are two groups of decisions, namely the aggregate level decisions and detailed level decisions, separated by the realization of uncertainties (Figure ~\ref{Fig-multi-stage procedure}).

Initially, a long-term staffing and scheduling plan is made, represented by the blue node in Figure ~\ref{Fig-multi-stage procedure}, followed by a sequence of adjustments. Within each stage $h$, there are two groups of decisions. Before the realization of uncertainties, aggregate level decisions for each stage $h$ aim to reduce the total cost by right-staffing (${n}_{h}$) through balancing under- and over-staffing costs ($\bar{n}_{h}$, $\tilde{n}_{h}$). On the detailed level, nurse scheduling decisions ${x}_{i,h}^{d,s}$ are more demand-adaptive, thus made after the uncertainties are revealed. Multi-stage stochastic programming is suitable for this problem because it considers the feature of making decisions sequentially \citep{birge2011introduction}.

\subsubsection{Discrete Scenarios}
\label{Assumption-discrete scenarios}
The uncertainty space can be captured by a finite set of discrete scenarios and the realizations of each stage are independent.

Now we introduce uncertainty realizations and scenarios of this problem. A scenario tree is adopted to describe the evolution process, assuming that the uncertainty space can be captured by a finite set of discrete scenarios denoted by $\Xi$. Uncertainty realization vector for each stage is represented by $ {\boldsymbol{\xi}}^{h}\in \Xi^{h}$, $ h\in \mathcal{H}$. 
An uncertainty realization is a tree branch. A sequence of uncertainty realizations until stage $h$, i.e., $({\boldsymbol{\xi}}^{1}, {\boldsymbol{\xi}}^{2},..., {\boldsymbol{\xi}}^{h})$, constitutes an $h$-stage scenario $\boldsymbol{\omega}_{h}\in\mathcal{W}_{h}$, corresponding to a tree node. A dummy 0-stage scenario $\boldsymbol{\omega}_{0}$ is assumed and no uncertainty is revealed at this stage, underlying its probability to be 1. Refer to the lower left side of Figure ~\ref{Fig-multi-stage procedure} for an illustrative example with two possible realizations for each stage. In this example, node $A$ is a scenario with uncertain realization $ {\xi}_{1}^{1}$ in first stage. According to Assumption ~\ref{Assumption-discrete scenarios}, we have $\mathcal{W}_{h}\equiv \Xi_{1}\times\Xi_{2}\times ... \times \Xi_{h} $ and $\mathcal{W}_{h}\equiv\mathcal{W}_{h-1}\times\Xi_{h}$, $h \in\mathcal{H}$, where $\mathcal{W}_{h}$ denotes the set of $h$-stage scenarios. Besides, let $\mathbf{a}(\boldsymbol{\omega}_{h})$ denote the unique parent node of scenario $\boldsymbol{\omega}_{h}$ and $\boldsymbol{\mathcal{G}}(\boldsymbol{\omega}_{h})$ represent the set of child scenarios in the subsequent stage $h+1$. Further, let $\boldsymbol{p}(\boldsymbol{\omega}_{h})$ be the probability of the $h$-stage scenario $\boldsymbol{\omega}_{h}$ and $\Bar{\boldsymbol{p}}(\boldsymbol{\xi}_{h})$ be the probability of the uncertainty realization $\boldsymbol{\xi}_{h}$. According to the independence in Assumption ~\ref{Assumption-discrete scenarios}, we have $\boldsymbol{p}(\boldsymbol{\omega}_{h}) = \prod_{h^{'}}{\Bar{\boldsymbol{p}}(\boldsymbol{\xi}_{h^{'}})}$ and $\sum_{\boldsymbol{\omega}_{h} \in \mathcal{W}_{h}}{\boldsymbol{p}(\boldsymbol{\omega}_{h})=1}$. 

\subsubsection{Scenario Tree}
\label{Assumption-known scenario tree}
The scenario tree is fully known for the entire planning horizon.

We assume the overall scenario tree can be estimated. The scope of this paper is to capture the evolving uncertainties in the decision process. With the help of data-driven optimization, this study can be further improved with the gradually known information of uncertainties, which will be a consideration for our future study.

%\vspace{0.6cm}

The multi-stage stochastic programming model is formulated as shown in Model 3 in Appendix \ref{app-model3}.

The objective function ~\ref{model3-stoch-obj} minimizes the total cost in expectation across the planning horizon, consisting of the cost of the initial stage and the recourse for each subsequent stage. Constraint ~\ref{model3-stoch-stageh-obj} defines the initial staffing and scheduling plan, which is a set of constraints from the deterministic model.  Function ~\ref{model3-stoch-obj2} gives the recourse function for stage $h$, consisting of the cost of over-staffing, under-staffing, the penalty for adjustments from the initial schedule, the penalty for supply-demand mismatch, and the penalty for soft rule violations. 

To model the aggregate level decisions for each stage, we introduce Constraints ~\ref{model3-stoch-flow conservation}-~\ref{model3-stoch-num of nurses2}. Constraint ~\ref{model3-stoch-flow conservation} enforces the conservation rule of the number of nurses. Constraint ~\ref{model3-stoch-num of nurses1} assumes that the transfer of nurses to other wards is satisfied first.  Constraint ~\ref{model3-stoch-num of nurses2} enforces the capacity rule, that is the total number of nurses with assignments at stage $h$ should not exceed the number of nurses and this is also the linking constraint of aggregate level decisions and detailed level decisions.  Then, similar to the deterministic model, a set of scheduling Constraints ~\ref{model3-stoch-h1} - \ref{model3-violations} are given for the detailed-level decisions. The domains of the variables are defined in Constraints ~\ref{model3-stoch-variable domain n} - \ref{model3-variable domain_p3}.

Note that in objective function ~\ref{model3-stoch-obj2}  non-linearity is introduced by the terms $|x_{i}^{d,s}- {x}_{i,h}^{d,s}(\boldsymbol{\omega}_{h})|$. The sign constraints are satisfied here because this term is positive for a minimization problem \citep{Cornell} and can be easily linearized as shown in Model 4 in Appendix \ref{app-linear-stoch}.

%% file: 04_sol_approach.tex
\section{Solution approach}
\label{Section-solution method}
In this section, we introduce the solution approaches. First, dependent on the structural property of the model, an equivalent MILP can be built, which reduces the computational complexity and then solved by a commercial solver (such as Gurobi or CPLEX). Second, to obtain near-optimal solutions within reasonable computational times, a generative AI-guided approach is given. 
\subsection{Structural property}
\label{Section-solution method-MILP}

\begin{theorem}
\label{proposition-block separable}
The problem formulated by Equation \ref{model4-stoch-linear-obj} has the block-separable recourse structure. The proof is given in Appendix \ref{app-proof-block-sep}.
\end{theorem}

% \proof{Proof of Proposition 1.} 
%  The proof is given in the e-companion (EC.6).
%  \endproof

According to Proposition ~\ref{proposition-block separable}, the following effect can be observed. Assuming independent uncertain realizations between stages, the under- and over-staffing variables at stage $h$ only depend on the number of nurses. In this way, the same ${n}_{h}$ is optimal for any realizations of $\boldsymbol{\xi}_{h}$. Then, the detailed level decision on nurse scheduling is left here to be optimized, which now depends separately on each stage's staffing decision. That is to say, the multi-stage program has a similar structure to a two-stage program. Refer to the following.

% \begin{remark}
\textbf{Remark 1. }
\label{Remark 1}
\ref{model4-stoch-linear-obj} is equivalent to a two-stage stochastic program with recourse. Specifically, an extensive form of the aggregate level problems constitutes the first-stage model. The second stage is composed of the detailed level recourse functions for the entire scenario tree. The equivalent model (Model 5) is shown in Appendix \ref{app-model5}. The proof is presented in Appendix \ref{app-proof-remark-1}.

In Model 5, the objective function in Equation ~\ref{model5-2stage-obj} is composed of two parts. First is the total expected cost of the aggregate level decision-making, i.e., the cost of under- and over-staffing, weighted by the probabilities of each scenario. Second is the summation of recourse (detailed level) of all stages. $\boldsymbol{\mathcal{G}}(\boldsymbol{\omega}_{h})$ represents the set of descent scenarios in the subsequent stage $h+1$.

\subsection{Generative method}
\label{Section-solution method-GFN}
Solving the above model to optimality is still challenging for the following reasons:  First, the size of the model is extremely large for realistic problem instances due to the scenario tree size. Second, the problem is further complicated by a very large number of integer variables. Hence, direct implementation to tackle this problem with commercial solvers typically does not yield high efficiency. A naive alternative is to replace all random variables with their expected values as proxies. Yet this may lead to discrepancies between the expected value and the actual value of different scenarios. Consequently, the expected value solutions may not be satisfactory when measured with stochastic settings. Nonetheless, expected values are still a "nice" estimator in expectation and are important for efficiency. In this section, our goal is to learn a policy function that can be trained using the expected values, while still being robust to different scenarios in stochastic settings. 

A common approach is to find a single optimal solution that minimizes the total expected cost with respect to all scenarios. In contrast, for our generative approach,  a set of diverse candidate solutions is constructed that have low costs according to the objective function. To be specific, we aim to learn a generative model that has higher probabilities to obtain low-cost solutions. 

The sections \ref{Section-flow networks} to \ref{Section-Generative nurse scheduling approach} below describe the following:
\begin{enumerate}
    \item Generative Flow Network (GFN) notations
    \item Training Objective
    \item Proposed Generative Nurse Scheduling (GNS) approach
\end{enumerate}

% In the following Section ~\ref{Section-flow networks}, notations on flow networks are given. Then, the training objective is introduced in Section ~\ref{Section-training objective}. Finally, our proposed solution approach is explained in Section ~\ref{Section-Generative nurse scheduling approach}.

\subsubsection{Generative flow networks (GFN)}
\label{Section-flow networks}
GFN is a novel generative modeling framework that trains models to sample diverse solutions proportional to the reward \cite{Bengio2021a,bengio2023gflownet}. GFN is particularly good at sampling compositional, discrete objects like graphs or sets. These objects are generated by a sequence of actions $a$, one part after the other, and ends at a terminal state $x$. The most important feature that distinguishes GFN from classic reinforcement learning (RL) methods is that instead of trying to find reward-maximizing solutions, GFN learns a policy $\pi(a|s)$ to sample action $a$ from state $s$ and the probability of visiting a terminal state $x$, i.e., $P(x)$, is proportional to its reward $R(x)$.

GFN is based on the flow of the directed acyclic graph (DAG), which is composed of a set of nodes $s\in\mathcal{S}$ and a set of edges $a\in\mathcal{A}$ that links the nodes. Particularly, there is only one source node 0 and multiple sink nodes $x\in\hat{\mathcal{S}} $ and $\hat{\mathcal{S}}\subset \mathcal{S}$ be the subset of sink nodes.  Further, $\mathcal{A}(s)$ is the subset of all edges that link to all successors from node $s$. $T(s,a)=s^{'}$ represents that start from node $s$ and take edge $a$ will lead to node $s^{'}$. A path  ${\tau}\in\mathcal{P}$ consists of a sequence of nodes $s$. All feasible paths must start from the root node, and can all be assigned some amount of flow. We differentiated partial paths that end with intermediate nodes and complete paths that end with sink nodes $x$. $\mathcal{P}_{x}$ denote the subset of paths that ends at sink node $x$. Note that GFN is based on a non-bijective graph. In other words, there could be several paths going through the same node, and also, except for the root, one node may have different predecessors, i.e., $\left|\left\{(s,a)|T(s,a)=s^{'}\right\}\right|\geq 1$. Let $F(s)$ be the total flow going through node $s$. $R(x)>0$ is the outflow of node $x$. $F(s,a)$ is the flow coming from node $s$ by edge $a$. Interestingly, the reader may notice that a GFN is analogous to a Markov decision process (MDP) where the root and sink nodes correspond to the initial and terminal state, while an edge corresponds to an action and a path corresponds to a sequence of actions. The forward probability $P_F$ corresponds to the policy function, while the out-flow is analogous to the reward function. 
 %(refer to Table in the e-companion (EC.3) to this paper.).

% \begin{figure}
%     \centering
%     \caption{An example of flow network with one source node $s_{0}$ and two sink nodes $s_{3}$ and $s_{6}$.}
%     \includegraphics[trim = 0 0 200 -70, scale=0.1,clip]{Figure/fig_flownet.jpg}  
%     \label{fig-flownet}
% \end{figure}

For a flow network like this, the flow conservation condition for each node $s\in\mathcal{S}$ must be satisfied. This is ensured by the following,
\begin{equation}
    \label{eq-gfn-inout1}
    \sum\limits_{s^{'},a:T(s^{'},a)=s}{F(s^{'},a)} = R(s) + \sum\limits_{a\in\mathcal{A}(s)}{F(s,a)}, \quad \forall s\in\mathcal{S}
\end{equation}
assuming that $R(s)=0$ for interior nodes, and $A(s)=\emptyset$ for sink nodes.

\

Then, for this flow network, the following proposition is proven by \cite{Bengio2021a} (proof presented in Appendix \ref{app-proof-prop-2})
\begin{theorem}
    \label{proposition-bybengio}
    If $P(s)$ denotes the probability of state $s$ going from the initial state. Let $\pi$ be the policy generating paths from the initial state by sampling actions $a\in\mathcal{A}(s)$ with 
    \begin{equation}
        \label{eq-proposition-begio1}
        \pi(a|s)=\frac{F(s,a)}{F(s)}.
    \end{equation}
    Then this results in, 
    \begin{equation}
        \label{eq-proposition-begio2}
        P(x)=\frac{R(x)}{\sum_{x^{'}\in\hat{\mathcal{S}}}{R(x^{'})}}, \quad \forall x\in \hat{\mathcal{S}}
    \end{equation}
    That is to say, the probability of visiting a terminal state $x$ ($P(x)$) is proportional to its reward $R(x)$.
\end{theorem}

% \proof{Proof of Proposition 2}
%     \label{proof-gfn-proportional}
%     The proof is provided in the e-companion (EC.8).
% \endproof

So according to Proposition ~\ref{proposition-bybengio}, the policy function $\pi$ can map the reward $R(x)$ to a generative model that generates paths that end with $x$ with a probability that is proportional to its reward. 

\subsubsection{Training objective}
\label{Section-training objective}
In this section, trajectory balance loss \citep{malkin2022trajectory} is used as the training objective function for the above generative model.  First, the notation of consistent flow is introduced. Then, the definition of trajectory balance loss is given.
% Then, for each path $\tau$, we have
% \begin{equation}
%     \label{eq-gfn-1}
%     P(\tau) = \frac{F(\tau)}{Z}, \quad \forall \tau \in \mathcal{P}
% \end{equation}
% \begin{equation}
%     \label{eq-gfn-2}
%     Z = \sum\limits_{\tau\in\mathcal{P}}{F(\tau)}
% \end{equation}
% where $F(\tau)$ is the flow of this path $\tau$ and $P(\tau)$ is the probability distribution of this path $\tau$. ~\ref[eq-gfn-2]{Eq.(\ref{eq-gfn-2})} expresses the flow conservation that the total inflow from the root node $s_0$ into the network equals the sum of all possible flows by the paths coming out of this node $s_0$. Let $x$ represent all complete trajectories that end at sink nodes $x$.  If we take a look at a particular path $x$ and let $R(x)$ be the flow of sink node $x$, we have

% \begin{equation}
%     \label{eq-gfn-6}
%     R(x) = \sum\limits_{\tau\in\mathcal{P}_{x}}{F(\tau)}.
% \end{equation}
% According to ~\ref[eq-gfn-6]{Eq.(\ref{eq-gfn-6})} and ~\ref[eq-gfn-1]{Eq.(\ref{eq-gfn-1})}, 
% \begin{equation}
%     \label{eq-gfn-7}
%     P(x)= \frac{R(x)}{Z} = \frac{\sum\limits_{\tau\in\mathcal{P}_{x}}{F(\tau)}}{Z}. 
% \end{equation}
Assume that apart from the source node, a path $\tau$ visits several nodes $s\in\mathcal{S}(\tau)$.  Let $\tau(s^{'})$ denote a partial path that ends at node $s^{'}$ based on path $\tau$. That is to say, path $\tau$ and $\tau(s^{'})$ share all nodes from the source until node $s^{'}$. Let $e_s$ denote the only one predecessor of node $s$ in path $\tau$, then
\begin{equation}
    \label{eq-gfn-3}
    P(\tau) = \prod\limits_{s\in\mathcal{S}(\tau)}{P_{F}\left[\tau(s)|\tau(e_s)\right]}, \quad \forall \tau\in\mathcal{P}
\end{equation}
and $P_{F}$ can be a conditional probability formulation that we generate paths starting from the source node. Similarly, there exists a backward probability $P_{B}$ factorizing the path probability conditioned on a sink node,
\begin{equation}
    \label{eq-gfn-4}
    P(\tau)=\prod\limits_{s\in\mathcal{S}(\tau)}{P_{B}\left[\tau(e_{s})|\tau(s)\right]}, \quad\forall \tau \in\mathcal{P}_{x}
\end{equation}

A consistent flow indicates that the flow in the forward direction should equal that of the backward direction \citep{bengio2023gflownet}. According to the flow conservation in the network, the forward direction flow is equal to the flow in the backward direction. Then this can be translated into the following conditions,

\begin{equation}
\begin{aligned}
    \label{eq-gfn-5}
    Z\prod\limits_{s\in\mathcal{S}(\tau)}{P_{F}\left[\tau(s)|\tau(e_s)\right]} &= R(x)\prod\limits_{s\in\mathcal{S}(\tau)}{P_{B}\left[\tau(e_{s})|\tau(s)\right]},
    \\ 
    \quad\forall \tau \in\mathcal{P}_{x}
\end{aligned}
\end{equation}

If we compute the squared difference between the logarithms of the LHS and RHS for Eq.~\ref{eq-gfn-5}, 
\begin{equation}
    \label{eq-gfn-6}
    \begin{aligned}
        \mathcal{L}(\tau;\boldsymbol{\theta}) 
        =\log R(x) + 
        \sum\limits_{s\in\mathcal{S}(\tau)}{P_{B}\left[\tau(e_{s})|\tau(s);\boldsymbol{\theta}\right]}\\ - 
        \sum\limits_{s\in\mathcal{S}(\tau)}{P_{F}\left[\tau(s)|\tau(e_s);\boldsymbol{\theta}\right]} 
    \end{aligned}
\end{equation}
Then in optimal, $\mathcal{L}(\tau,\boldsymbol{\theta})$ is equal to $\log Z$, and is the same for all paths within the flow network. Thus, the training objective is to minimize the loss variance over all paths $\tau$,
\begin{equation}
    \label{eq-gfn-7}
    \mathcal{L}_{obj}(\tau,\boldsymbol{\theta}) = {\left( \mathcal{L}(\tau,\boldsymbol{\theta}) -\mathbb{E}_{\tau}[\mathcal{L}(\tau,\boldsymbol{\theta})] \right)}^{2}, \quad \forall \tau \in\mathcal{P}
\end{equation}
Note that in the above training, only forward and backward probabilities are needed to learn. And after training, we can directly generate solutions with policy $P_{F}$, which is proportional to the reward.

\subsubsection{Generative nurse scheduling approach (GNS)}
\label{Section-Generative nurse scheduling approach}
\begin{figure*}
    \centering
    \caption{Generative nurse scheduling approach}
    \includegraphics[width=0.9\linewidth]{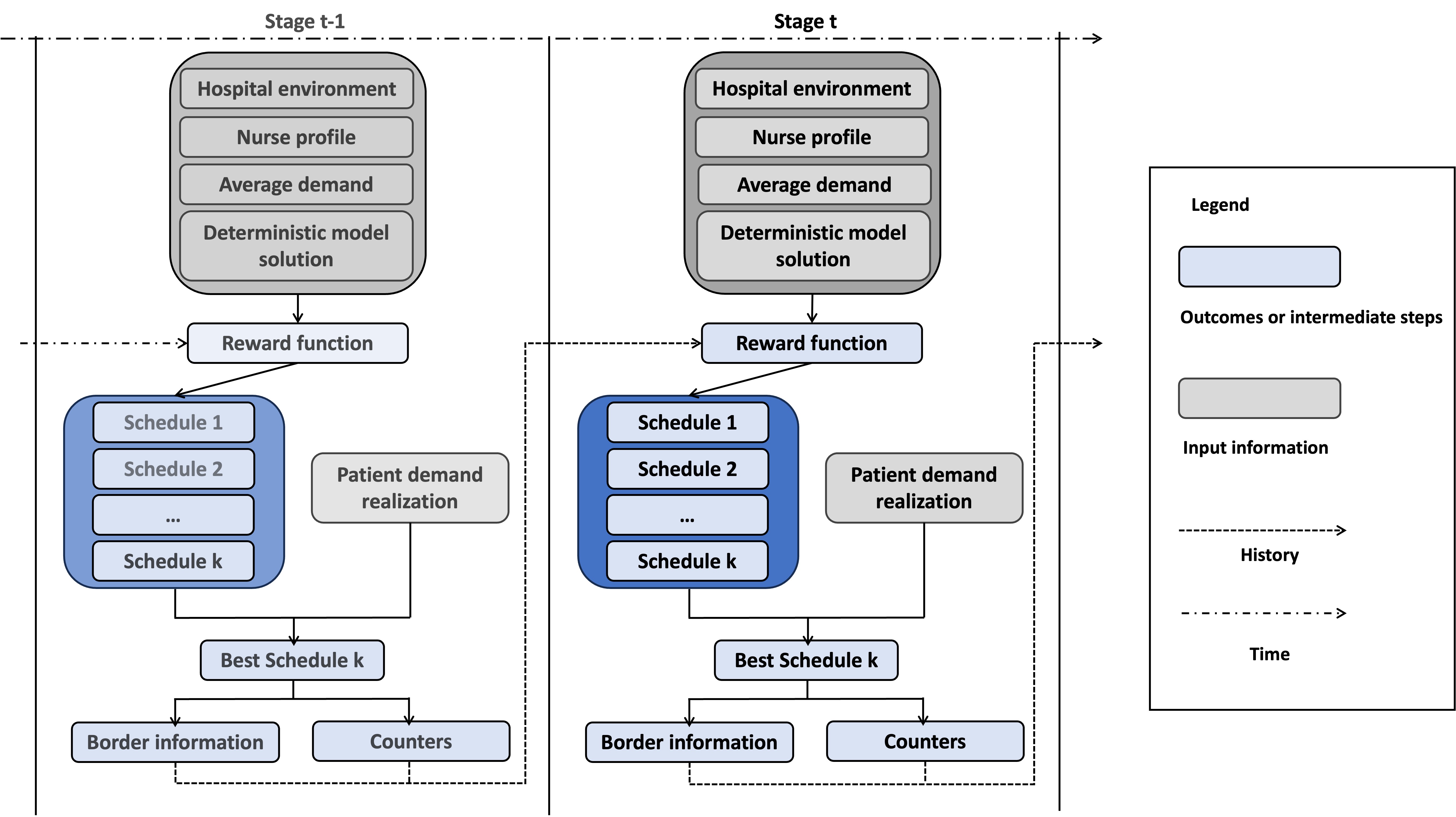}
    \label{fig-rolling}
\end{figure*}

% For the proposed nurse scheduling problem with evolving uncertainties, exploration is especially useful since the entire roster can be generated from stage to stage, considering all new information. 
Based on the concept of GFN discussed above, this section presents our generative nurse scheduling (GNS) approach in detail. This approach has two main features. First, utilizing a generative model, the nurse schedule is obtained via a sequence of actions, i.e., making an assignment for a nurse. Each generated schedule corresponds to a terminal state in the GFN as introduced in \ref{Section-flow networks}. Second, we obtain a good set of candidate staffing and scheduling decisions relying on the expected value of stochastic demand before the realization of uncertainties at each stage. The framework of GNS is shown in Figure ~\ref{fig-rolling}. At each stage, with the hospital environment, nurse profile, average demand and deterministic model solution as the input information, we take Objective ~\ref{model-proxy} as the reward function to train the model. The deterministic \ref{model2-obj} solution is utilized as an initial solution and several adjustments are made to it. Since the deterministic model is not robust to future information considering only a single scenario, two enhancements are introduced here. First, embedded with GFN, using the trained model, we can generate diverse solutions and pick the top $k$ schedules as candidates, which improves the robustness of the method. Second, after the realization of uncertain patient demand, we undergo an evaluation process that assesses all good-performing candidate solutions with real patient demand. Since this evaluation process computes the cost of each schedule without optimization, it is expected to be efficient computationally. Then, the border information and corresponding counters are passed to the next stage. 
In this setting, the reward function $R(x)$ is as follows,

\begin{align}
    \label{model-proxy}
    R\left(x| \boldsymbol{\chi}_{0}, \bar{\boldsymbol{D}}, \rho\right) &= M - 
    \frac{\left(U_{2}\left( \boldsymbol{\chi}_{0} \right) + \sum_{h\in\mathcal{H}}{\Upsilon_{h}\left(  {\boldsymbol{\chi}}^{h-1}, \bar{\boldsymbol{D}}\right)}\right)}{\rho}
\end{align}
where we map the primal minimization problem into maximization since we use the reward here. $\bar{\boldsymbol{D}}$ denotes the expected demand and the reward can be regarded as the simplified version of Function \ref{model5-2stage-obj} with a single scenario. $\rho \in \mathbb{R}^{+}$ represents the temperature, which improves the efficiency of exploration and avoids to stuck into the sub-optimal. This is useful since the number of schedules with higher reward values is much less than that of the lower-reward ones. Without using the temperature, it is likely that a large set of samplings is required to find a good schedule. The temperature $\rho$ enables us to make a trade-off between diversity and shifting sampling distributions towards better solutions \citep{zhang2023robust}.

\begin{figure*}
    \centering
    \caption{Actions sampling and schedule. For the current step, the action NurseB-Tuesday-P2 is drawn from the action space and added to the current partial schedule (Suppose P2 is one of the PM shifts, A1 is one of the AM shifts, N1 is one of the N shifts). }
    \includegraphics[trim = 0 0 0 -70, scale=0.09,clip]{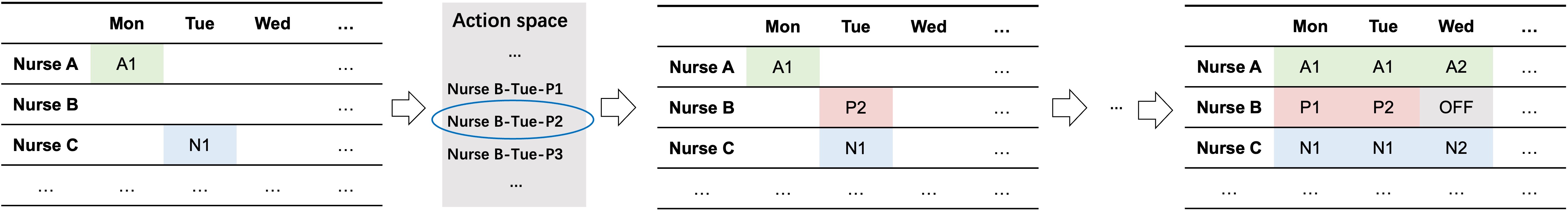} 
    \label{fig-gfn-nsp}
\end{figure*}

Specifically, the main part of GNS follows a step-by-step process within one episode (\ref{alg:gns_method}). First, each episode starts with an empty schedule as the initial state. Then, at each step, according to the policy function, an action is sampled and added to the schedule. The action here is the assignment (as shown in Figure ~\ref{fig-gfn-nsp}). Note that the set of valid actions at each step is defined and updated based on all constraints in Model 5 in Appendix \ref{app-model5} and all actions that will lead to infeasible partial schedules are masked out. The sampling is repeated until the predefined termination and each episode returns a schedule and with its corresponding reward. For every $t$ episode, based on the training objective loss function ~\ref{eq-gfn-7}, the gradient descent is conducted, and the hyper-parameters are updated. 

\begin{figure*}

\begin{algorithm}[H]
\label{alg:gns_method}
\caption{GNS method}
\KwIn{\textnormal{Hospital environment configuration, Nurse profile, Average demand, Deterministic Model 2 and its initial solution $x_0$, num of epochs, num of actions $T$, update\_freq}}
\KwOut{\textnormal{A policy function}}
Initialization: an empty schedule as \texttt{InitialSchedule}\;
Define Policy to predict $P_F$, $P_B$ and $\log Z$ with neural networks\;
episode $\leftarrow 0$, stop $\leftarrow$ false\;
\While{episode $<$ num of epochs}{
    \tcc*{Step 1: Start with initial schedule}
    State $\leftarrow$ InitialSchedule\; 
    $P_F, P_B \leftarrow \text{Policy}$(State)\; \tcc*{Step 2: Predict using Policy}
    \For{$i \leftarrow 1$ \KwTo $T$}{
        cat $\leftarrow \text{Categorical(logits=}P_F)$\; \tcc*{Sample actions based on $P_F$}
        action $\leftarrow \text{cat.sample()}$\;
        NewState $\leftarrow$ State + [keys[action]]\; \tcc*{Step 4: Transition to new state}
        Total$_{P_F}$ $\leftarrow$ Total$_{P_F}$ + cat.log\_prob(action)\;
        \If{$i = T$}{
            Reward $\leftarrow$ reward(NewState)\; \tcc*{Terminal state, compute reward}
        }
        $P_F, P_B \leftarrow \text{Policy}$(NewState)\; \tcc*{Recompute policy}
        State $\leftarrow$ NewState\;
    }
    Compute losses based on Equations 58 and 61\;
    \If{episode mod update\_freq == 0}{
        Update neural network parameters\; \tcc*{Gradient step}
    }
    episode $\leftarrow$ episode + 1\;
}
\end{algorithm}
\end{figure*}

%% file: 05_experiments.tex
\section{Numerical experiments and case study}
\label{Section-empirical study}
Now, we present numerical results based on anonymized data from a surgical ward of the Tan Tock Seng Hospital (TTSH) in Singapore. First, the data overview and current practice are introduced in Section ~\ref{Section-Data overview} and Section ~\ref{Section-Current roster overview}. Then, the computational performance is analyzed in Section ~\ref{Section-Computational performance}.  In Section ~\ref{Results of reward values for GNS}, the results of reward values for the GNS are given. Major results of real-world instances are given in Section ~\ref{Section-Results of real-world cases}. Then, the effect of the bounded flexibility is tested in Section ~\ref{Section-case-Trade-off}. To verify the significance of considering uncertainties, the value of the stochastic solution is analyzed in Section ~\ref{Section-VSS}. All algorithms are implemented in Python on a personal computer with a 13th Gen Intel(R) processor and 16 GB RAM. 

\subsection{Data overview}
\label{Section-Data overview}
TTSH comprises clinical, allied health departments, and specialist centers, and employs thousands of healthcare staff. It ranks second in Singapore's large acute care general hospital with more than 1,500 beds and the busiest trauma center in the country \citep{TTSHwiki}. The 5-month dataset we obtained is from a particular surgical ward with 60 nurses, which contains the predicted daily patient demand, nurses' preferred working times and requests, as well as the prevailing scheduling rules for that planning horizon. 
%In real-world cases, the nurses can be transferred between different wards.

\subsection{Current scheduling process overview}
\label{Section-Current roster overview}
There are 12 types of different shifts as shown in Figure ~\ref{fig-shift timings and time slots in the studied ward}. The lengths of shifts vary between 6 to 12 hours. Five belong to the AM slot with red blocks, and six belong to the PM slot with green blocks. There is only one shift in the night slot (N), which is colored blue. The studied ward schedules their nurses by the following four steps. First, forecasting. They predict demand for each time slot and day. Second, nurse staffing. The number of nurses required to cover the demand is determined. Third, collecting nurses' requests. The information on nurses' requests is collected. Some nurses' particular shifts are approved. Fourth, rostering. Nurses are assigned to day and shifts.  Nurses are rostered based on a pre-determined set of hard and soft rules. In the rostering process, the following challenges are faced by the nursing manager:

\begin{figure}
    \centering
    \caption{Shifts and time slots.}
    \includegraphics[width=0.9\linewidth]{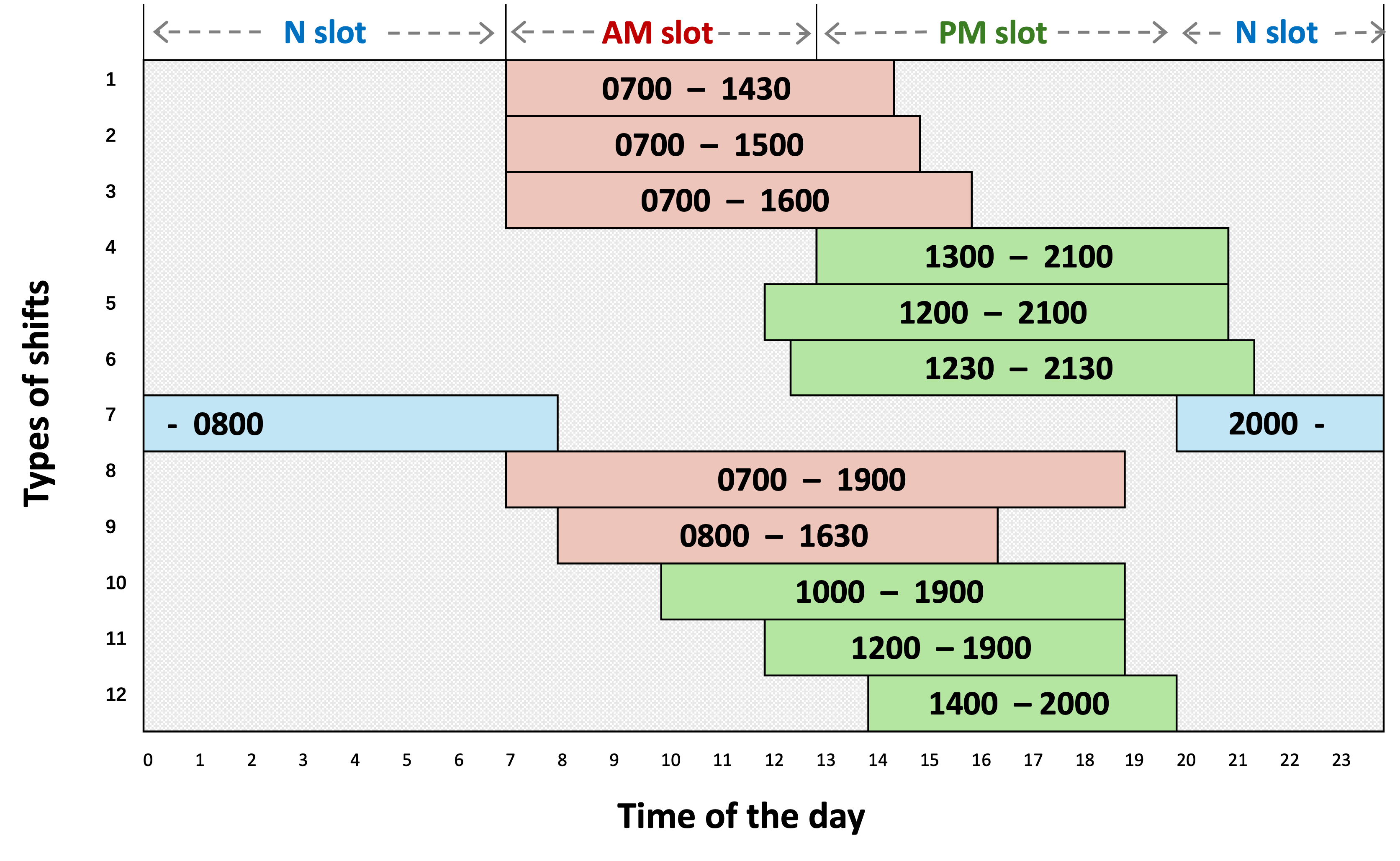}
    \label{fig-shift timings and time slots in the studied ward}
\end{figure}

First, based on a traditional manpower scheduling tool that does not handle the full-fledged set of nurse preferences, rosters are manually revised via the experience of nurse managers. Consequently, the rostering process is quite labor-intensive and error-prone.  

Second, there are mismatches of demand and supply. Especially, the mismatches of different time slots are quite different. For example, on particular days, there can be over-supplied nurses on the AM slot and under-supplied nurses on the PM slot. A clear objective function is required to describe and minimize the so-called imbalanced situation.
Third, to enhance the welfare of nurses, many hard and soft rules are included. Some of them are rarely seen in the cases of hospitals and even contradict each other. The logic behind each rule is required to be studied.

\subsection{Computational performance}
\label{Section-Computational performance}
\begin{figure}
    \centering
    \caption{Scenario tree in use.}
    \includegraphics[trim = 0 0 500 -100, width=0.9 \linewidth,clip]{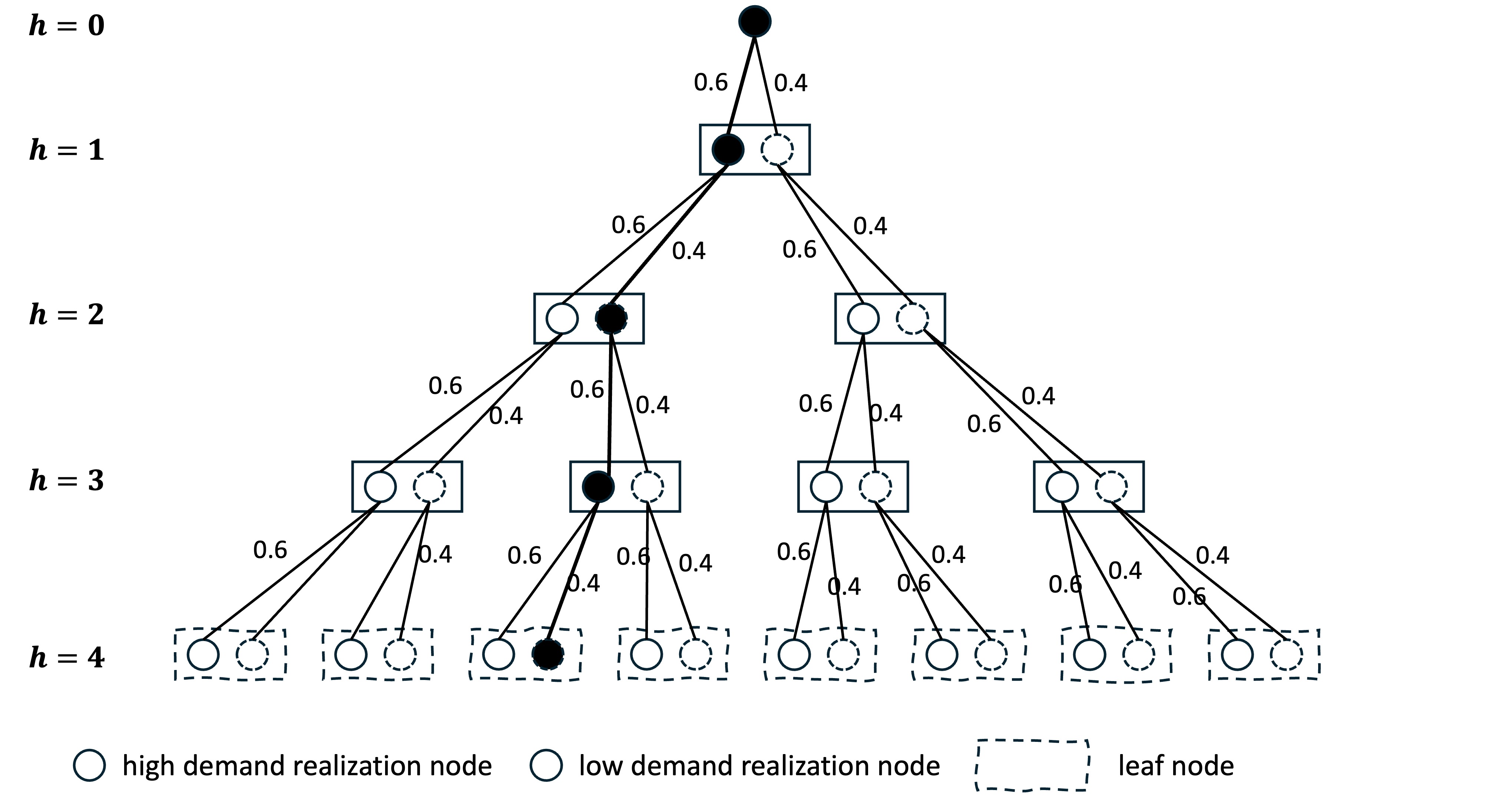}
    \label{fig-scenario tree in use}
\end{figure}

In this section, we present the computational performances of our proposed approach. The main results are shown in Table ~\ref{table-computational performance}. The first column of the table lists various cases. To facilitate comparisons of various settings and problem scales, we generate 28 cases, each representing a different combination of the size of nurses' sets $|I|$ and the demand scale $\eta$, denoted by $|I|$-$\eta$. For example, Case 20-0.5 is the case with an overall 20 nurses and a demand scale of 0.5. Here is how we simulate each case. The demand scale is an input parameter that describes different levels of patient demand. For each case, the initial demand $\mathbf{D}^{0}$ is set to be 5, 6 and 4 for the AM, PM and N slots respectively each day. Then the scaled demand is calculated as  $\eta\mathbf{D}^{0}$. Two possible demand realizations, high (+4 demand each slot) and low (-2 demand each slot) are assumed with probabilities $p_{h}=0.6$ and $p_{l}=0.4$ for all stages as described in Figure ~\ref{fig-scenario tree in use}. Then, the expected demand $\mathbf{D}^{E}$ for each stage is $p_{h}(\eta\mathbf{D}^{0}+4\mathbf{e})+p_{l}(\eta\mathbf{D}^{0}-2\mathbf{e})$ rounded up to an integer, where $\mathbf{e}$ represents the vector of 1s. For example, in case 10-0.5, the expected demand is $0.6*(2.5+4,3+4,2+4)+0.4*(2.5-2,3-2,2-2)\approx(6,5,4)$. Next, for each stage, the high-demand and low-demand realization is randomly generated by  $rand(\mathbf{D}^{E}+4\mathbf{e}, 2)$ and $rand(\mathbf{D}^{E}-2\mathbf{e}, 2)$, where $rand(a,b)$ denotes a normal distribution with mean $a$ and standard deviation $b$. To rule out the effect of stochasticity, each case is composed of 5 instances by randomly generating the demand 5 times. Thus there are a total of 140 instances. 

The second to fifth columns show the number of variables, and solution time (in seconds) for the deterministic program (DP) and stochastic program (SP), respectively. From columns 6 to 8, the objective values (total cost) of DP, SP,  and GNS are given. Finally, the optimality gaps of SP and GNS are provided in the last two columns. For SP, Table ~\ref{table-computational performance} lists the average solution time and objective value of the 5 instances. An upper time limit of 7,200 CPU seconds is set and the current best solution found will be returned if hitting the time limit. For the GNS method, the overall training episode is set to be 50,000 in each case.

Note that the DP objective value measures the total cost with a single demand scenario whereas the SP and GNS measure demand uncertainty and the corresponding recourse, and hence their values show a significant difference. And since the DP can always be solved to optimality, the gap of DP is not listed. The gap for SP is the gap returned from the commercial solver within the time limit (i.e. the branch and cut upper and lower bounds). On the other hand, the gap for GNS is the relative difference between the GNS objective value and the corresponding SP objective value, which can be negative, since the SP may not be solved to optimality within its time limit. 

\begin{figure*}
        \centering
    \caption{Results of the reward distributions}
 \includegraphics[trim = 0 0 0 0, width=0.9 \linewidth,clip]{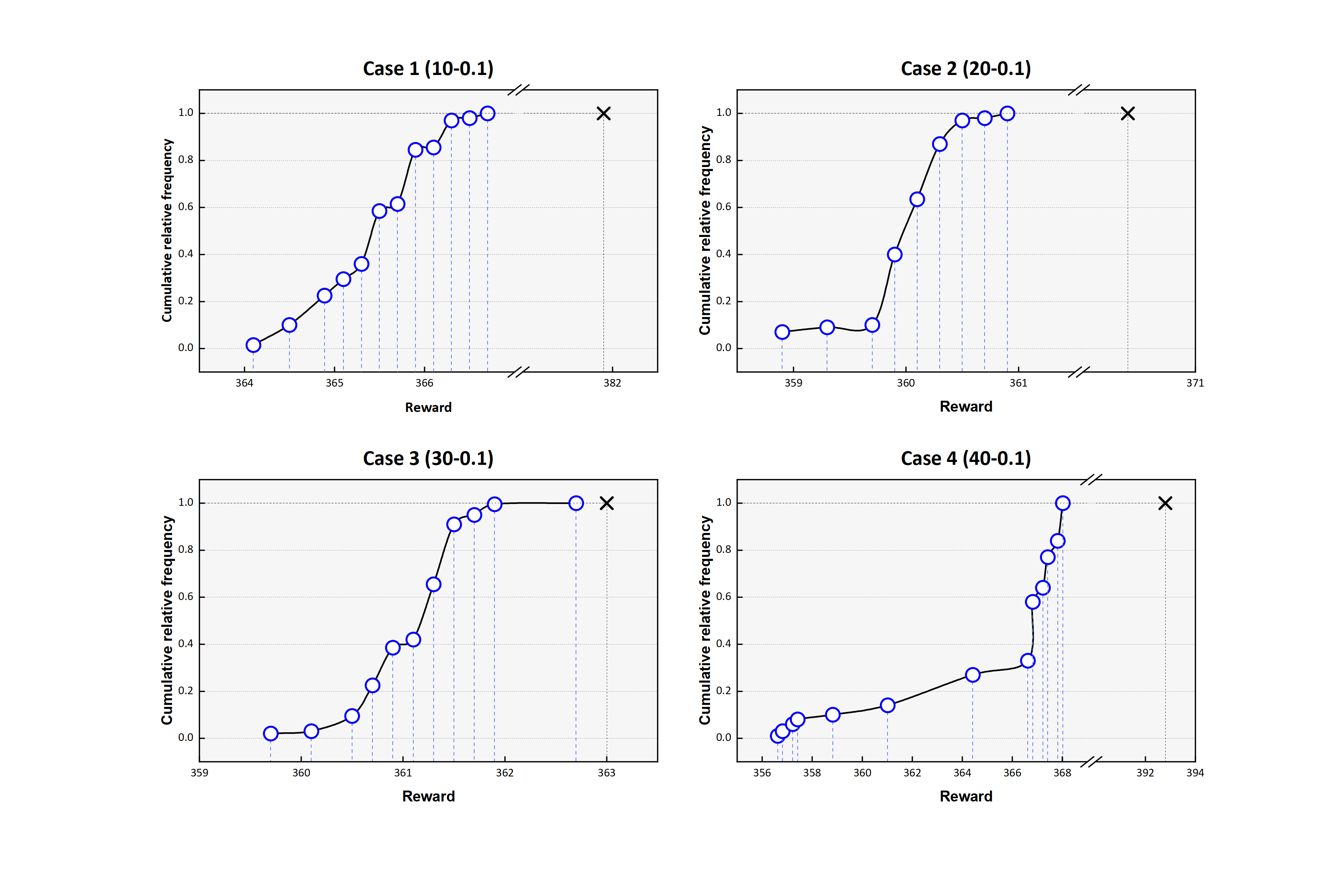}\label{Fig-prob_rw}
\end{figure*}

As shown in Table ~\ref{table-computational performance}, due to the multi-stage structure, the number of variables of SP is much larger than DP for all cases. The SP also used much more CPU time (ranging from 83.08 sec. to an upper limit of 7200 sec) against different cases. The solution time for SP hits the two-hour time limit with a total number of 20 nurses and the DP can be solved to optimal within seconds. We can also observe that the gap of SP ranges from 0\% to as extremely high as 519.33\%. Note that even though the average SP solution time for cases 20-0.5, 20-1 and 20-1.5 does not exceed the upper limit, there exist some instances that were solved to feasibility instead of optimality among the 5 instances, due to the time limit. Hence, the gap is not 0 for those 3 cases, but rather, they are less than 2.52\%, which is acceptable for real-world applications. Since the key idea of the SP is that we can optimize solutions all at once assuming that we can somehow know or predict the scenario tree, the solution time limit can be further lifted and the commercial solver can still be utilized to generate the optima or near-optima based on the hospital's requirement. And here we provide an alternative heuristic, which is using the generative model to progressively generate the schedule. It can be observed from the table that GNS shows great potential if we increase the problem scale. For some cases with a supply of 40 nurses, GNS outperforms the solver of SP, with a negative gap as high as -9.28\%. Note that for cases 40-0.01, 40-0.05, 40-0.1, 40-0.25, the gaps of SP are significantly large, mainly because those are the cases with large supply relative to the demand which are challenging for the solver within the time limit; whereas the GNS approach shows some promising trend in handling large instances well.

\begin{table*}
\caption{Performance of different solution approaches.}
\centering
        {\def\arraystretch{1.4}  %change the number for increasing or decreasing the spacing.
            %\tiny
            \scriptsize %uncomment this for changing the font size
\begin{tabular*}{0.97\textwidth}{@{\extracolsep{\fill}}llllllllll}
\hline
\hline
  Case & \multicolumn{2}{c}{Num variables} & \multicolumn{2}{c}{Solution time (sec)} & \multicolumn{3}{c}{Objective value (cost)} & \multicolumn{2}{c}{GAP (\%)} \\ \cline{2-3} \cline{4-5}\cline{6-8} \cline{9-10} 
$|\mathcal{I}|$-$\eta$ & DP & SP & DP & SP & DP & SP & GNS& SP & GNS \\ \hline
10-0.01 & 5,131  & 75,197&  0.21   &813.63 &204.00&5744.73&6012.94  & 0.00& 4.67 \\
10-0.05& 5,131 & 75,197 & 0.26  &1399.73&204.00 &6538.14&6895.16 & 0.00 &5.46 \\
10-0.1& 5,131 & 75,197 &0.46  & 1356.74& 297.00 & 6569.34& 8043.29 &0.00 & 22.44\\
10-0.25 & 5,131&  75,197
      & 0.37 &  1311.90 &1761.00 &8796.72 & 8846.21 & 0.00  & 0.56 \\
10-0.5 &5,131  & 75,197
      &0.21  &2002.84 &5121.00 &  14122.06&  14679.04 & 0.00  & 3.94 \\
10-1 & 5,131  & 75,197
      &0.23 &400.11  & 11001.00 &  25786.55&26013.27 & 0.00 &0.88 \\
10-1.5&5,131 & 75,197
      &0.28  &83.08  &17721.00   &39350.02 & 40157.38 & 0.00  &2.05  \\
20-0.01 &  10,177 & 149,007&  0.64   &5,840.83 &405.00&2,318.63&2712.17  & 0.00& 16.97 \\
20-0.05& 10,177 & 149,007 & 0.58  &4,076.73&405.00 &2,441.87& 2978.62& 0.00 & 21.98\\
20-0.1 & 10,177& 149,007 & 0.62 & 4213.58 & 405.00 & 2160.33 & 3028.13 & 0.00 &40.17\\
20-0.25 & 10,177&  149,007
      & 0.66 &  4,499.29 &405.00 &3,471.59 & 3,565.24 & 0.00  & 2.70 \\
20-0.5 &10,177  & 149,007
      & 1.18 & 6,626.23&423.00 &  5,165.01&  5,240.01 & 2.52  & 1.45 \\
20-1 & 10,177  & 149,007
      &1.27 & 6,106.64 & 4,332.00 &  12,618.40& 12,942.25 & 0.05 & 2.57 \\
20-1.5& 10,177 & 149,007
      & 1.46 &  7,186.62&  11,052.00 &25,070.39 & 25,128.06 & 0.11  &  0.23\\
30-0.01 & 15,277  & 222,817
      &0.71 &  7,200.00&  504.00 & 497.86& 653.16 & 56.60  & 31.19 \\
30-0.05 &15,277  &222,817
      &1.02 &  7,200.00&504.00 &831.06  &835.29  & 14.90  & 0.51\\
30-0.1 & 15,277  & 222,817 & 1.12& 7200.00 & 504.00 & 552.25 &836.78
& 57.67 & 51.52 \\
30-0.25 & 15,277 & 222,817
      &1.16 &   7,200.00&  504.00 &704.78& 743.04 & 85.22  & 5.43 \\
30-0.5 &15,277 &  222,817 
      & 1.37 & 7,200.00&  504.00 &2,324.18& 2,348.04 & 53.52 & 1.03  \\
30-1 & 15,277&  222,817 
      &2.23 &  7,200.00&  513.00& 5,261.15& 5,692.05 & 15.91 &8.19   \\
30-1.5 &15,277 &  222,817 
      & 2.16 &7,200.00 &4,020.00 &12,323.73 & 12,375.10 & 5.09  &0.42 \\
40-0.01 &20,377 &  296,627 
      & 1.02 &  7,200.00&519.00 & 357.18 & 486.04  & 488.49 &   36.08\\
40-0.05 & 20,377&  296,627 
      &1.39   &7,200.00& 519.00  &378.97&  526.17&507.69 & 38.84 \\
40-0.1 & 20,377 &  296,627 & 1.08 & 7200.00 & 519.00 & 567.38 & 734.76 & 481.16 & 29.50 \\
40-0.25 &20,377  & 296,627 
      &1.08  &7,200.00 &519.00 &445.60 &  628.24&519.33 &  40.99\\
40-0.5 & 20,377 & 296,627 
      &1.19  & 7,200.00& 519.00& 670.66 &  619.40&307.64 & -7.64  \\
40-1 & 20,377 & 296,627
      &1.19   &7,200.00 &519.00 &1,816.00 &  1,745.52&116.50 & -3.88 \\
40-1.5 &20,377 & 296,627 
      &3.53  &7,200.00 & 537.00&6,315.12 &  5,729.30&32.35 & -9.28\\
\hline
\hline
\end{tabular*}
\label{table-computational performance}
}
\end{table*}

\subsection{Results of reward values for GNS}
\label{Results of reward values for GNS}
GNS is capable of generating solutions that are proportional to the reward values. After training the model, we output the final 500 schedules generated by GNS and analyze their reward values' distribution. Here we present the cumulative relative frequency \footnote{Cumulative relative frequency is the ratio of the number of solutions whose reward values equal or exceed the respective values on the x-axis.} distributions obtained in various illustrative cases (Figure \ref{Fig-prob_rw}). We map backward the single reward value obtained from the stochastic programming(SP) model and the result is presented by the X node in the figure. We can make two observations from the figure. First, the high reward values constitute a majority of all obtained reward values. For example, in case 1 (i.e. 10-0.1), over 70\% of reward values are higher than 366. This shows that the policy tends to have a higher probability of generating a decent solution. Second, for all those four cases, the slope for the cumulative  reward distribution increases with higher reward values, indicating that the training policy is capable of generating higer reward value solutions with higher probabilities.

\subsection{Results of real-world instances}
\label{Section-Results of real-world cases} 
\begin{figure}
        \centering
    \caption{Results of daily demand coverage gap by time slot.}
 \includegraphics[trim = 0 0 0 40, scale=1,clip, width= \linewidth]{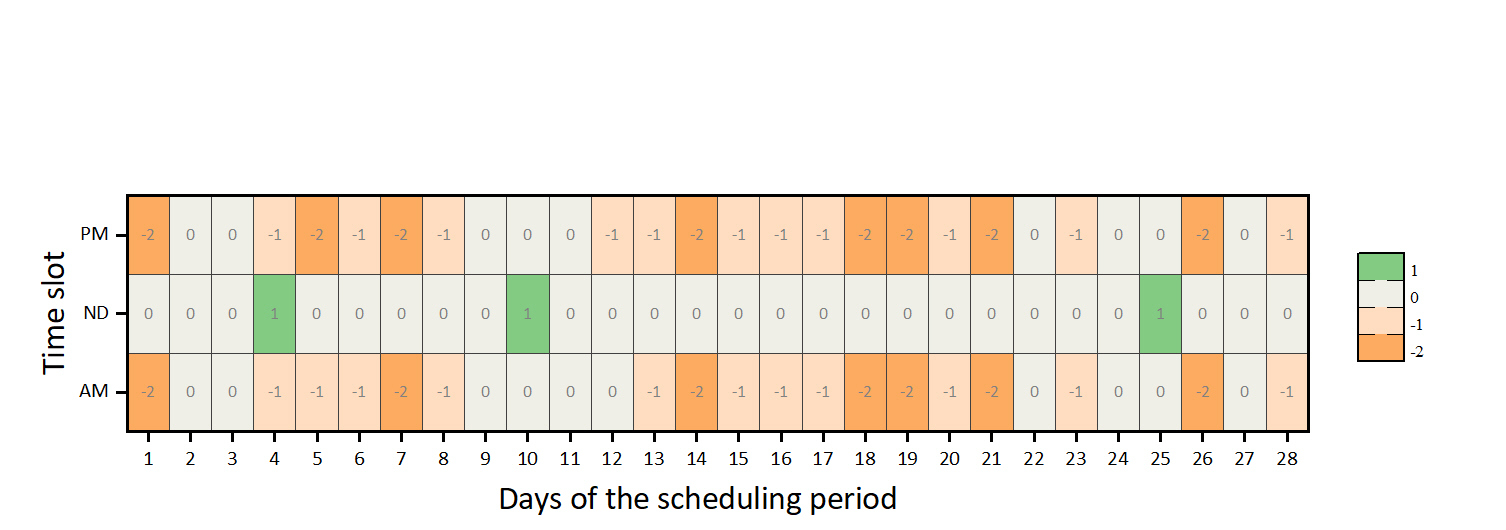}\label{Fig-daily demand coverage gap}
\end{figure}

\begin{figure}
        \centering
    \caption{Results of daily demand coverage gap by skill.}
 \includegraphics[trim = 0 250 350 180, width=0.9 \linewidth,clip]{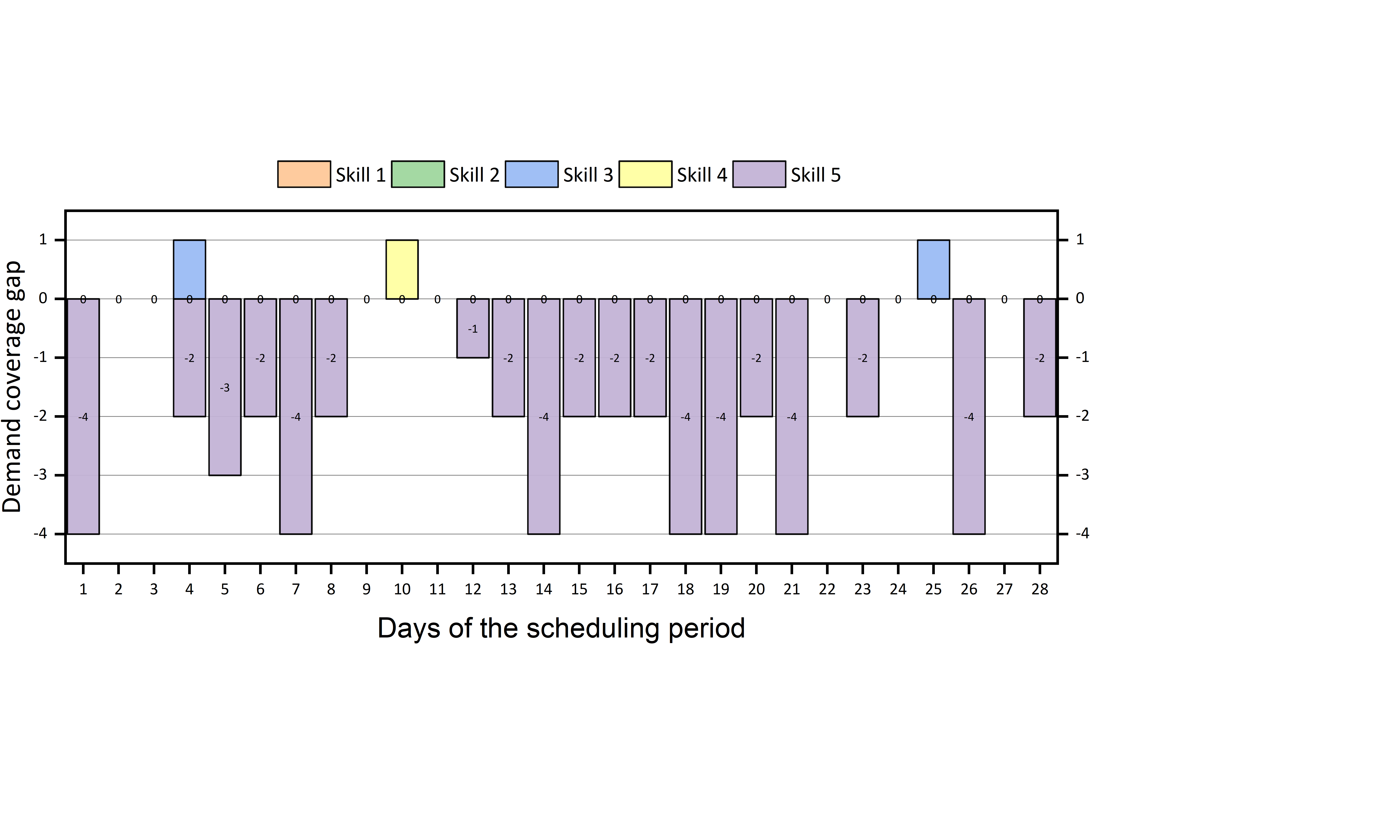}\label{Fig-coverage_skill}
\end{figure}

We now present the results of real-world cases. Figure \ref{Fig-daily demand coverage gap}
shows the gap for daily demand coverage in each time slot. The x-axis represents different days of the scheduling period and the y-axis is different time slots within a day. For each block, there is a number representing the difference between demand and supply (demand coverage gap). We observe that the gap has been restricted in $[-2,1]$, where we have a maximum of 2 understaffing and 1 overstaffing among 47 total mismatches. On average, this equals a 1.68 absolute gap on average per day. To further strengthen the demand coverage constraint, the demand coverage by skill is added to the model based on the real-world requirements of the studied hospital. Nurses are divided into five categories, from Skill 1 to Skill 5 (from the most senior to the most junior) and the cost of understaffing for senior nurses is much higher than that of junior nurses. To show the result of demand coverage by skills, Figure \ref{Fig-coverage_skill} is given. The stacked column graph shows the sum of the demand coverage gap for all time slots by different skills. It can be observed from the figure that there is no understaffing of senior nurses, with only 1 overstaffing on Skill 3 and Skill 4. The understaffing cases are all due to the Skill 5 group.

Taking a closer look at the results of the daily demand coverage gap (in Figure \ref{Fig-daily demand coverage gap}), we observe that the absolute gaps with each time slot are very close on a single day (the maximum absolute difference is 1). The reason is that we further impose the "balance" constraint here. In the real-world application, the study hospital faces challenges in the too-busy morning and too-free afternoon. The benefit can be achieved by introducing "balance" to represent the total absolute gap difference between each time slot. That is to say, the daily demand coverage gap can be further improved by a more balanced distribution of the mismatch cases. The proposed model outperforms the current practice as it includes the requirements more precisely. 

\subsection{Results of bounded flexibility}
\label{Section-case-Trade-off}
\begin{figure}
        \centering
    \caption{Results of bounded flexibility}
 \includegraphics[trim = 0 50 0 100, width=0.9 \linewidth,clip]{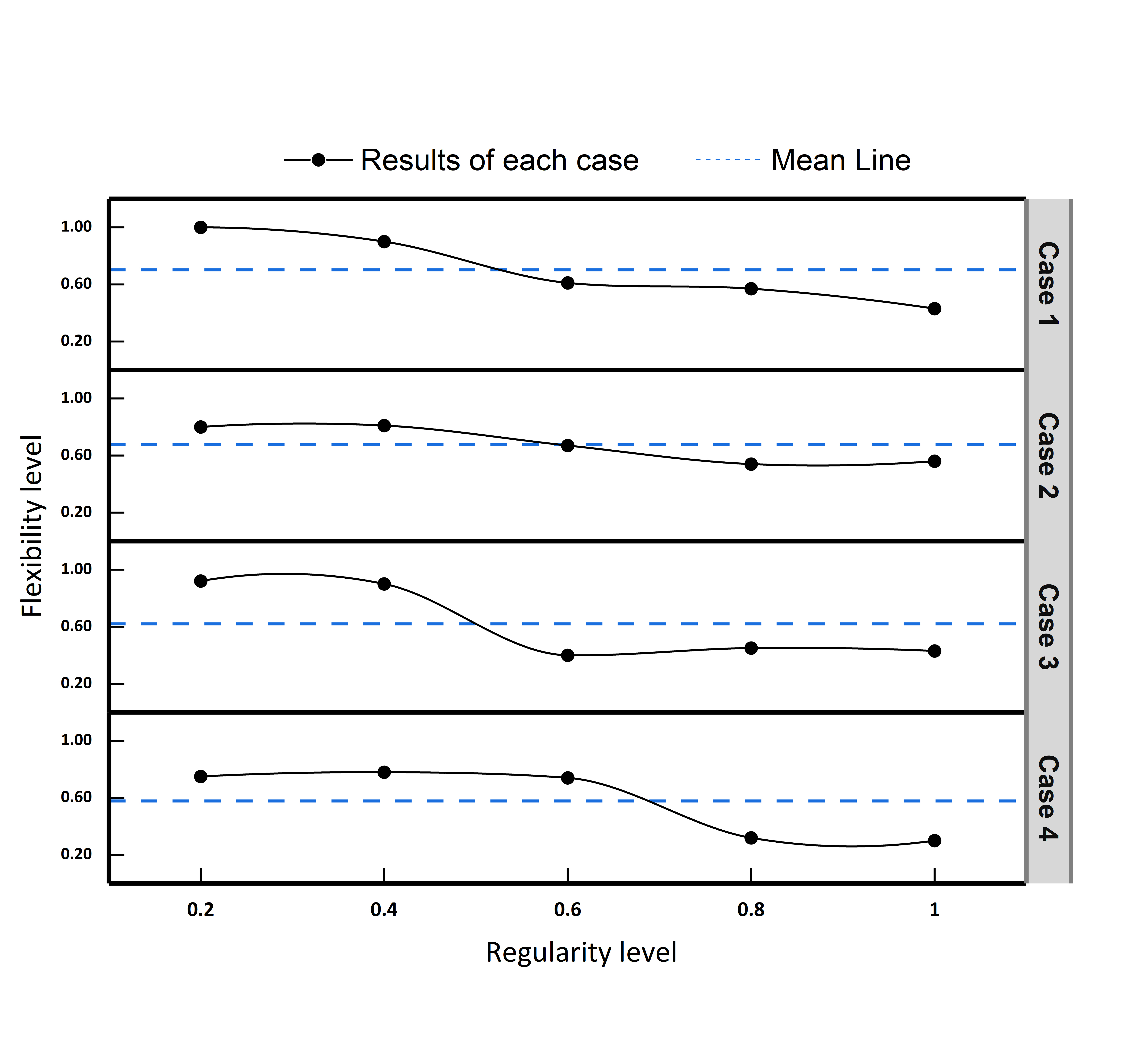}
 \label{Fig-flexibility-regularity}
\end{figure}

In this section, we present the effect of bounded flexibility. Against the work policy, we show the results of flexibility level under different regularity level settings. In the following experiments, flexibility level refers to the satisfaction rate of all nurse requests and regularity level is the ratio of nurses with $p_{1}$ or $p_{2}$ work policy. Note that we assume that $p_{3}$ is the most irregular work policy since the nurse can have shifts fall into any time slot in the scheduling horizon. Four cases of different combinations of total nurse number and demand scale are tested and the results are shown in Figure \ref{Fig-flexibility-regularity}. The blue lines present the mean flexibility level for each case. It can be observed from the figure that with the increase in regularity, the flexibility takes on a decreasing trend, from around 0.8 to 0.3. Also, the results reveal the average flexibility levels under different regularity settings are 0.71, 0.68, 0.61 and 0.60, respectively. 

The following conclusions can be drawn. First, increasing the amount of nurses with $p_{1}$ and $p_{2}$ work policies reduces the satisfaction rate of nurses' requests. In other words, restricting the component of nurses with $p_{3}$ work policy may sacrifice request satisfaction. The main reason is that nurses themselves may submit requests that already go against the regularity of the schedule. $p_{3}$ nurses can be assigned any time slot during the scheduling period, which creates flexibility. Results also show that a small reduction in the regularity level from 0.6 to 0.4 could remarkably bring higher flexibility, which could possibly be achieved by hiring more nurses with $p_{3}$ policy and across specialties. 

Further, the average flexibility level remains stable in different cases. The reason is that we aim to minimize the request rejection rate in the objective function, so the model avoids the situation with very low flexibility.

\subsection{Value of stochastic solution}
\label{Section-VSS}
Most studies associated with the nurse scheduling problem rely on either deterministic or two-stage stochastic models, in order to avoid the modeling challenges and solving complexity. To further quantify the value of using the multi-stage stochastic program, two benchmarks are introduced. First, the value of the stochastic solution is measured by the comparison with a deterministic model based on the expected value of uncertainties. Second, the value of sequential decision-making is verified by the comparison with a two-stage stochastic programming.

The Value of Stochastic Solution (VSS) \citep{birge2011introduction} is measured by the difference between the Expected Value Problem (EVP) and the multi-stage stochastic programming problem. First, the EVP is defined as,

\begin{equation}
\label{eq-VSS-EVP}
U_{EVP}=\min_{\chi}U(\chi,\textbf{D}^{E})
\end{equation}
where $\chi$ is decision variable vector and $\textbf{D}^{E}=\mathbb{E}[\textbf{D}]$  represents the expectation of the random patient demand $\textbf{D}$. With the above model solution, let $\bar{\chi}_{1}(\textbf{D}^{E})$ denote its optimal solution for the initial roster and all aggregate level decision variables in the subsequent stages. Now, we use VSS to measure the performance of an EV solution against various scenarios in the detailed level decision-making process, i.e., assigning nurses to day and shift. Also, the expected result of the EV solution (EEV) is,

\begin{multline}
	\label{eq-EEV}
	U_{EEV}
 =\min_{\chi_{2}} U_{0}\left( \bar{\chi}_{1}(\textbf{D}^{E}) \right)
  + C_{S}n_{0}
  \\ + 
\mathbb{E}_{\mathbf{D}} \Bigg [Q_{1}\left( \bar{\chi}_{1}(\textbf{D}^{E}),\mathbf{D}\right) + ... 
 \\ +  \mathbb{E}_{\mathbf{D}}\left[Q_{|\mathcal{H}|}\left( \bar{\chi}_{1}(\textbf{D}^{E}),\mathbf{D}\right)\right]...
    \Bigg]
\end{multline}

Let $\chi_{2}$ denote the detailed-level decision vector. The EEV can reflect how $\bar{\chi}_{1}(\textbf{D}^{E})$ perform under all possible random scenarios, considering the detailed-level decisions in optimal as functions of $\bar{\chi}_{1}(\textbf{D}^{E})$ and $\mathbf{D}$. Finally, the VSS is computed by,

\begin{equation}
\begin{gathered}
    \label{eq-VSS-final}
	VSS = U_{EEV} - U_{5}^{*}
\end{gathered}
\end{equation}
where $U_{5}^{*}$ denotes the optimal value for the proposed multi-stage stochastic programming problem. For the comparison with a two-stage stochastic programming model, we present a simplified model where all aggregate level decisions are made before the planning horizon considering the recourse function with 16 scenarios (if we take the example of a four-stage problem with 2 scenarios each stage). This underlies an assumption that the future demand realization of each stage is fully known by us and there are no sequential decisions since the demand is not sequentially revealed.

\begin{table}
\caption{Value of stochastic solutions.}
\centering
        {\def\arraystretch{1.5}  %change the number for increasing or decreasing the spacing.
            \scriptsize %uncomment this for changing the font size
\begin{tabular*}{0.48\textwidth}{@{\extracolsep{\fill}}llllll}
\hline
\multirow{2}{*}{\textbf{Case}: $|\mathcal{I}|$-$\eta$} & \multirow{2}{*}{PP} & \multirow{2}{*}{EEV} & \multirow{2}{*}{TP} & \multicolumn{2}{c}{VSS} \\ \cline{5-6} 
& & & & EEV-PP & TP-PP\\ \hline
% 5 - 0.01&8739.53&9114.82&9070.99&375.28&331.46\\ 
% 5 - 0.05&9221.46&9544.41&9447.17&322.95&225.71\\ 
% 5 - 0.25&14289.33&14671.00&14450.31&381.67&160.98\\
% 5 - 0.5&21270.09&21397.33&21289.13&127.24&19.04\\ 
% 5 - 1&33188.85&33254.99&33188.85&66.14&0.00\\ 
% 5 - 1.5&46675.41&46676.91&46675.41&1.50&0.00\\ 
10-0.01& 5744.73& 6862.75&6452.73 &1118.02 & 708.00\\
10-0.05&6538.14 &7680.32 & 7265.62&1142.18 &727.47 \\ 
10-0.1& 6569.34&7210.38 &7130.88 & 641.04 & 561.54\\
10-0.25& 8796.72&9315.34 & 9245.07&518.62 & 448.35\\ 
10-0.5&14122.06 & 14552.72& 14335.64& 430.66&213.58 \\ 
10-1&25786.55 &26000.76
 &25788.86 &214.21 & 2.30\\
10-1.5&39350.02 & 39450.31 & 39350.02&100.29& 0.00\\
20-0.01&2318.63&5909.32&4581.96&3590.69&2263.33\\
20-0.05&2441.87&6093.66&4590.74&3651.78&2148.87\\
20-0.1 & 2160.33 & 5076.00
& 4071.39&2915.67&1911.06\\
20-0.25&3471.59&5862.05&5021.48&2390.46&1549.89\\
20-0.5&5165.01&6429.67&6219.29&1264.66&1054.28\\
20-1&12618.40&13153.98&13340.71&535.58&722.31\\
20-1.5&25070.39&25378.70&25374.21&308.30&303.82\\
30-0.01 &497.86&4478.32&3335.04&3980.46&2837.17\\
30-0.05 &831.06&5563.86&4004.31&4732.81&3173.26\\
30-0.1 &552.25&4738.78&3825.50&4186.52&3273.25 \\
30-0.25 &704.78&4769.48&3485.10&4064.70&2780.32\\
30-0.5&2324.18&5721.29&4616.77&3397.10&2292.59\\
30-1&5261.15&6713.32&9130.16&1452.17&3869.01\\
30-1.5&12323.73&12515.78&12933.70&192.04&609.97\\
40-0.01&357.18&4549.57&4806.74&4192.38&4449.56\\
40-0.05&378.97&5299.56&4085.70&4920.59&3706.73\\
40-0.1&567.38&5218.55&4799.88&4651.18&4232.50\\
40-0.25&445.60&4776.13&3569.04&4330.52&3123.43\\
40-0.5&670.66&4748.70&4059.60&4078.04&3388.93\\
40-1&1816.00&5040.49&4187.15&3224.48&2371.15\\
40-1.5 &6315.12&7213.96&8545.62&898.84&2230.50\\
\hline
\end{tabular*}
\label{table-VSS}
}
\end{table}

The results are listed in Table \ref{table-VSS}, where the first column shows different cases and the second to fourth columns are the objective values of the multi-stage stochastic model (PP, we simplify it to be the primal problem), EEV and two-stage stochastic model (TP), respectively. The last two columns list the VSS. The parameter setting is identical to the above sections. For each case, we run 5 instances and return the average values. It can be observed from the table that all VSSs are positive (the last two columns), meaning that a nurse manager can yield hidden savings considering demand uncertainty. Also, under the same level of demand scale, we yield larger VSSs with more total number of nurses. For example, the figure rises from 3590.69 and 2263.33 in case 20-0.01 to 4192.38 and 4449.56 in case 40-0.01.  In other words, compared with the EEV, the PP and the TP can both generate robust solutions to different scenarios and lower total costs. Positive values of VSSs can reflect that with demand uncertainty, the decision-making process is able to balance all possible realizations to save the total costs. While for the EVP, only single one scenario is included. Another observation is that PP outperforms TP, with more VSS generated for each case. In the PP, the latter-stage decisions can be optimized after getting to know the uncertainty under the former-stage decisions, and the consideration of evolving uncertainty enables the model with more flexibility. Consequently, the higher VSSs can be seen in column EEV-PP than in TP-PP. On the opposite way, under the same total number of nurses, the VSSs decrease as demand rises. In the cases with 20 nurses, the values drop from 3590.69 and 2263.33 to 308.30 and 303.82 as the demand scale goes from 0.01 to 1.5. To summarize, the cost of ignoring the demand uncertainties is notably much higher than incorporating the stochastic and fluctuating patient demand into planning and operational decisions in nurse management. Besides, the cost of a simplified two-stage model is significantly greater than that of a model with demand-evolving features. So it is essential for nurse managers to incorporate this effect for long-term benefits. Together, results suggest that it is of equal importance to consider demand uncertainty and its evolving features.

%% file: 06_conclusion.tex
\section{Conclusions}
\label{Section-conclusions}
This study incorporates bounded flexibility and demand uncertainty in nurse scheduling, motivated by a real-world problem faced by a Singapore hospital. A deterministic nurse staffing and scheduling model is presented that coordinates various rostering rules and bounded flexibility to help accommodate nurse-specific needs. Then, given evolving uncertainties in patient demand, a multi-stage stochastic programming model is formulated with the aim of minimizing the total system costs. The key to solving this problem is the utilization of the block-separable structure that remarkably reduces the complexity of the model. Another contribution is the design of a state-of-the-art generative scheduling approach, by which we can generate schedules that are proportional to their costs. This is especially effective and has managerial relevance where the diversity of schedules is expected. From the practical perspective, the real-world case study generates the following insights. 
First, the computation of the value of stochastic solutions shows that the model can achieve remarkable savings considering the evolving uncertainties. Second, the proposed formulation also makes use of the clear distinction between hard rules that should be strictly respected and flexible soft rules that support bounded flexibility. A small reduction in the regularity level could remarkably bring higher nurse flexibility. Third, the daily demand coverage gap can be further improved by a more balanced distribution of the mismatch cases. The proposed model outperforms the current practice as it includes the requirements more precisely. 

This work is aimed at motivating more research interest in stochastic workforce optimization. A possible extension is the consideration of supply-side uncertainties. Nurses' absence can have a great impact on the operations and sometimes it is challenging for the nurse managers to adjust and fix the schedules with all the requirements. Empirical studies of how nurse-specific features (e.g., designation and experience) will be useful for estimating nurses' absence or requests. A better understanding of the factors that influence a nurse’s supply at the individual level can help derive best practices so as to improve operational performance at the system level. Furthermore, since after using the block-separable structure, the problem becomes a MILP, it may further improve the computational performance with exact algorithms to speed up the MILP.

%% file: 07_appendix.tex
\clearpage

\appendix

\section{Abbreviations}
\label{app-abbrev}
\begin{table}[ht]
\caption{\textbf{Abbreviations}}
\centering
        {\def\arraystretch{1.5}  %change the number for increasing or decreasing the spacing.
            \scriptsize %uncomment this for changing the font size
\begin{tabular*}{0.4\textwidth}{@{\extracolsep{\fill}}ll}
\hline
\hline
\textbf{Abbreviation} & \textbf{Definition} \\
\hline
AM &  morning slot\\
PM &  afternoon slot\\
N &   night slot\\
NW &  nonworking day\\
W &  working day\\ 
\hline
\hline
\end{tabular*}
\label{table-abbreviations}
}
\end{table}

\section{Notations for deterministic model}
\label{app-notations-deterministic}
\setlength{\headheight}{45pt} 
\vspace{0.2cm}
% \begin{center}
% \centering
        {\def\arraystretch{1.5}  %change the number for increasing or decreasing the spacing.
            %\scriptsize %uncomment this for changing the font size
% \begin{tabular*}{1\textwidth}{@{\extracolsep{\fill}}ll}
\begin{tabular}{ m{0.2 \linewidth} m{0.75 \linewidth} }
\multicolumn{2}{l}{\textbf{Sets and indices}}\\
$\mathcal{I}$ & set of nurses, $i \in \mathcal{I}$ \\
$\mathcal{D}$ & set of days, $d \in \mathcal{D}$ \\
$\mathcal{S}$ & set of shifts, $s \in \mathcal{S}$ \\
$\mathcal{J}$ & set of time slots, $j \in \mathcal{J}=\{AM, PM, N\}$\\
$\mathcal{P}$ & set of work policy, $p \in \mathcal{P}$ \\
$\mathcal{K}$ & set of weeks, $k \in \mathcal{K}$ \\
$\mathcal{I}_{p}$ & subset of nurses with work policy $p$\\
$\mathcal{S}_{j}$ & subset of shifts belongs to slot $j$  \\
$\mathcal{S}_{i}$ & subset of preferred shifts of nurse $i$\\
% $\mathcal{S}_{w}$ & subset of shifts belongs to team $w$, i.e., flexible shift and base shift\\
$\mathcal{D}^{L}_{t}$ & subset of last $t$ days of the planning horizon\\
% $\mathcal{D}^{F}_{t}$ & subset of first $t$ days of the planning horizon\\
$\mathcal{D}_{k}$ & subset of days at week $k$ \\
% $\mathcal{D}_{W}$ & subset of weekend days (Saturday and Sunday)\\
$\mathcal{D}_{Sat}, \mathcal{D}_{Sun}$ & subset of Saturday, Sunday\\
\end{tabular}

\begin{tabular}{ m{0.2 \linewidth} m{0.75 \linewidth} }
\multicolumn{2}{l}{\textbf{Parameters}}\\
% $A^{max}, A^{min}$& maximum/minimum number of N shifts in the planning horizon\\ 
$N$ &total number of available nurses\\
$T^{min}_{i}$& minimum work hours for nurse $i$\\
$T^{max}_{i}$& maximum work hours for nurse $i$\\
$E_{i}$& minimum non-working days per week for nurse $i$\\
 % $WD$ & number of weekend days in the planning horizon policy  $p$ in the planning horizon \\
 $L_{s}$ & number of effective work hours for shift $s$ \\
 % $B_{i}^{p}$& 1 if the work policy of nurse $i$ is $p$, 0 otherwise\\ 
 % $BL_{i}^{d}$& 1 if nurse $i$ requests a nonworking day $d$, 0 otherwise\\ 
$R_{i}^{d,s}$& 1 if nurse $i$ requests to work shift $s$ on day $d$, 0 otherwise\\ 
$W_{p}$& number of time slots for work policy $p$ \\ 
% $N_{p}$ & maximum number of nights for nurse with work policy $p$ in the planning horizon\\
$F_{i}$ & maximum number of consecutive working days for nurse $i$\\
$V_{i}$ & maximum number of violations for nurse $i$\\
$M_{i}$ & maximum number of working weekend days for nurse $i$\\
$Q_{j}^{d}$& number of nurses needed to cover patient demand for slot $j$ on day $d$ \\
% $C_{R}$ & weight for the fairness objective\\
$C_{C}$ & penalty for per demand-supply mismatch\\
$C_{R}$ & penalty for per request rejection\\
$C_{S} $ & cost for scheduling one nurse in the planning horizon\\
$C_{m}$ & penalty for having in total $m$ violations \\
\end{tabular}

\begin{tabular}{ m{0.2 \linewidth} m{0.75 \linewidth} }
\multicolumn{2}{l}{\textbf{Variables}}\\
$\beta_{i}$ & 1 if nurse $i$ is scheduled in the planning horizon (0-1)\\
 $x_{i}^{d,s}$& if nurse $i$ is assigned to shift $s$ on day $d$ (0-1)  \\
 $y_{i}^{j}$ & if nurse $i$ is assigned to time slot $j$ in the planning horizon (0-1)\\
 $\bar{q}_{j}^{d}$ & under-supplied demand for slot $j$ on day $d$\\
 $\tilde{q}_{j}^{d}$ &over-supplied demand for slot $j$ on day $d$\\
 $v_{i}^{m}$ & 1 if nurse $i$ has a total number of $m$ violations in the planning horizon (0-1)\\
 $p_{1,i}$ & number of days exceeding the maximum working weekend days $M_{i}$ for nurse $i$\\
$p_{2,i}^{d}$ & 1 if nurse $i$ has an undesired pattern NW-AM from day $d$ (0-1)\\
% $p_{2,i}$ & number of nights over the maximum value $N_{p}^{max}$ for nurse $i$ with work policy $p$\\
% $p_{3,i}^{d}$ & 1 if appearance of undesired pattern ND-ND-W on three consecutive days of nurse $i$ from day $d$\\
% $p_{2,i}$ & number of days over the maximum consecutive working days $F_{i}$ of nurse $i$\\
$p_{3,i}^{d}$ & 1 if nurse $i$ has an undesired pattern N-N on weekends from day $d$ (0-1) \\
$p_{4,i}^{d}$ & 1 if nurse $i$ has an undesired pattern NW-W-NW  from day $d$ (0-1)\\
$Y_{i}^{d,s}$& auxiliary variable for the linearization of requests rejection (0-1)\\
\end{tabular}
\label{table:sets and indices}
}
% \end{center}

\clearpage

\input{07_app_model1}

\input{07_app_model2}

\input{07_app_model3}

\input{07_app_model4}

\input{07_app_model5}

-
\clearpage

% \section{Linearization of the deterministic model}
% \label{app-linear-deterministic}
% [Model 2]
% \setlength{\jot}{0.3pt}

% \begin{multline}
%     \label{model2-obj}
%     minimize \quad U_{2}\left( \boldsymbol{\chi} \right) 
%  =\overbrace{C_{S}\sum\limits_{i\in\mathcal{I}}{\beta_{i}}}^{Staffing} 
%  + \overbrace{C_{C}\sum\limits_{j\in\mathcal{J}}{\sum\limits_{d\in\mathcal{D}}{(\bar{q}_{j}^{d}+\tilde{q}_{j}^{d})}}}^{Coverage}
%  \\+ \overbrace{
%             C_{R}\sum\limits_{i\in\mathcal{I}}{\sum\limits_{d\in \mathcal{D}}{\sum\limits_{s\in\mathcal{S}_{i}}{R_{i}^{d,s}Y_{i}^{d,s}}}}}^{Requests} 
% + 
% \overbrace{\sum\limits_{i\in\mathcal{I}}{
% \sum\limits_{m=1}^{V_{i}}{C_{m}v_{i}^{m}}
% }}^{Violations}
% \end{multline}

% \begin{equation}
% \label{model2-linearizations}
% \begin{cases}
% Y_{i}^{d,s} \leq 1 - x_{i}^{d,s} & \\ 
% Y_{i}^{d,s} \leq \beta_{i} & \\ 
% Y_{i}^{d,s} \geq \beta_{i} - x_{i}^{d,s} & \quad \forall i\in\mathcal{I}, s\in \mathcal{S}_{i}, d \in \mathcal{D}\\
% Y_{i}^{d,s} \in \{0,1\} & 
% \end{cases}
% \end{equation}

% \begin{equation}
%     \label{model2-others}
%     \text{Constraints} \quad ~\ref{model1-num of nurses2}-~\ref{model1-variable domain_p3}
% \end{equation}

\section{Notations for stochastic model}
\label{app-notations-stoch}
\setlength{\headheight}{55pt} 

\begin{center}
\label{table:stochastic-sets and indices}
\centering
        {\def\arraystretch{1.5}  %change the number for increasing or decreasing the spacing.
            %\scriptsize %uncomment this for changing the font size
\begin{tabular}{ m{0.2 \linewidth} m{0.75 \linewidth} }
\multicolumn{2}{l}{\textbf{Sets and indices}}\\
$\mathcal{H}$ & set of stages, $h=0, 1, ...,   |\mathcal{H}|$ \\
$\mathcal{D}_{h}$ & subset of days belong to stage $h$ \\
$\mathcal{D}^{h,L}_{t}$ & subset of last $t$ days for stage $h$\\
% $\mathcal{D}^{h,F}_{t}$ & subset of first $t$ days for stage $h$\\
\end{tabular}

\begin{tabular}{ m{0.2 \linewidth} m{0.75 \linewidth} }
\multicolumn{2}{l}{\textbf{Parameters}}\\
$C_{U}$& cost of under-staffing, i.e., outsourcing nurses\\
$C_{O}$& cost of over-staffing, i.e., canceling some nurse-shift  \\
$C_{A}$& penalty for per schedule adjustment\\
$T_{i, h}^{min}$& minimum work hours for nurse $i$ at stage $h$\\
$T_{i, h}^{max}$& maximum work hours for nurse $i$ at stage $h$\\
$E_{i, h}$& minimum non-working days for nurse $i$  at stage $h$\\
$F_{i,h}$ & maximum number of consecutive working days for nurse $i$  at stage $h$\\
$V_{i, h}$ & maximum number of violations for nurse $i$  at stage $h$\\
$M_{i, h}$ & maximum number of working weekend days for nurse $i$  at stage $h$\\
$Q_{j, h}^{d}$& number of nurses needed to cover patient demand for slot $j$ on day $d$  at stage $h$\\
\end{tabular}

\begin{tabular}{ >{\centering\arraybackslash} m{0.2 \linewidth} m{0.75 \linewidth} }
\multicolumn{2}{l}{\textbf{Variables}}\\
Stage 0: & \\
$\boldsymbol{\chi}$ & initial staffing and  scheduling decisions, i.e., $\boldsymbol{\chi}=\left\{\boldsymbol{\beta},\mathbf{x,y,\bar{q},\tilde{q}}\right\}$ as specified in the deterministic model. \\
\end{tabular}

\begin{tabular}{ m{0.2 \linewidth} m{0.75 \linewidth} }
\multicolumn{2}{l}{Aggregate level:}\\
${n}_{h}(\boldsymbol{\omega}_{h})$ & number of nurses scheduled at stage $h$\\
$\bar{n}_{h}(\boldsymbol{\omega}_{h})$ & number of under-staffing nurses at stage $h$\\
$\tilde{n}_{h}(\boldsymbol{\omega}_{h})$ & number of over-staffing nurses at stage $h$\\
% $ {\alpha}_{i}^{h}$ & percent of request fulfilled of nurse $i$ on week $h$ \\
\end{tabular}

\begin{tabular}{ m{0.2 \linewidth} m{0.75 \linewidth} }
\multicolumn{2}{l}{Detailed level:}\\
 $ {x}_{i,h}^{d,s}(\boldsymbol{\omega}_{h})$& if nurse $i$ is assigned to shift $s$ on day $d$ at stage $h$ (0-1)\\
$ {\beta}_{i,h}(\boldsymbol{\omega}_{h})$ & if nurse $i$ is scheduled at stage $h$ (0-1)\\
$ \bar{q}_{j,h}^{d}(\boldsymbol{\omega}_{h})$ & number of under-supplied demand in slot $j$ on day $d$ at stage $h$\\
$ \tilde{q}_{j,h}^{d}(\boldsymbol{\omega}_{h})$ &number of over-supplied demand in slot $j$ on day $d$ at stage $h$\\
$v_{i,h}^{m}(\boldsymbol{\omega}_{h})$ & 1 if nurse $i$ has a total number of $m$ violations at stage $h$ \\
$p_{1,i}^{h}(\boldsymbol{\omega}_{h})$ & number of days over the maximum working weekend days $M_{i,h}$ for nurse $i$ at stage $h$ \\
$p_{2,i}^{d,h}(\boldsymbol{\omega}_{h})$ & 1 if nurse $i$ has undesired pattern NW-AM  from day $d$ at stage $h$ \\
$p_{3,i}^{d,h}(\boldsymbol{\omega}_{h})$ & 1 if nurse $i$ has an undesired pattern N-N on on weekends from day $d$ at stage $h$\\
$p_{4,i}^{h}(\boldsymbol{\omega}_{h})$ & 1 if nurse $i$ has an undesired pattern NW-W-NW from day $d$ at stage $h$\\
$U_{i,h}^{d,s}$& auxiliary variable for the linearization of adjustments (0-1)\\
% $p_{5,i}^{h}(\boldsymbol{\omega}_{h})$ & number of days less than the minimum number of nonworking weekend days $NW_{i}^{min}$ for nurse $i$\\
% $p_{6,i}^{d,h}(\boldsymbol{\omega}_{h})$ & 1 if the appearance of undesired pattern ND-ND  on weekends for nurse $i$ from day $d$  \\
% $p_{7,i}^{d,h}$ & 1 if appearance of undesired pattern NW-W-NW for nurse $i$ starting from day $d$ \\
\end{tabular}
}
\end{center}

\newpage

\section{Proof of Proposition \ref{proposition-block separable}}
\label{app-proof-block-sep}

 \label{proof-block separable}
    As described by \cite{birge2011introduction}, a multi-stage stochastic linear program has block-separable recourse if the following four conditions are satisfied:
    First, for all stages $h\in\mathcal{H}$ and scenarios $\omega_{h}$, the decisions $x_{h}(\omega_{h})$ can be divided into aggregate level decisions ${z}_{h}(\omega_{h})$ and detailed level decisions ${y}_{h}(\omega_{h})$. Second, a stage $h$ objective contribution is in the form of $c_{h}x_{h}(\omega_{h})=r_{h}{z}_{h}(\omega_{h}) + q_{h}{y}_{h}(\omega_{h})$. Third, the recourse matrix $W_{h}$ is block diagonal:
    \begin{equation}
        \label{eq-proof-block-1}
        W_{h} = 
        \begin{pmatrix}
            A_{h} & 0\\
            0 &  B_{h}
        \end{pmatrix}
    \end{equation}
    Fourth, for $h$-stage scenario $\omega_{h}$, the technology matrix $T(\omega_{h})$  and right-hand-side matrix $H(\omega_{h})$ satisfying:
    \begin{equation}
        \label{eq-proof-block-2}
        T(\omega_{h}) = \begin{pmatrix}
            R(\omega_{h}) & 0\\
            S(\omega_{h}) &  0
        \end{pmatrix} \quad and \quad H(\omega_{h}) = 
        \begin{pmatrix}
            b(\omega_{h})\\
            d(\omega_{h})
        \end{pmatrix}
    \end{equation}
    where 
    $W_{h}x(\omega_{h}) = H(\omega_{h}) -T(\omega_{h-1})x_{h-1}$.

In \ref{model4-stoch-linear-obj}, within each stage $h$, there are two levels of decisions, namely staff sizing from the aggregate level and nurse scheduling from the detailed level. Besides, for each stage, the recourse function ~\ref{model4-stoch-obj} is composed of the cost of aggregate level decisions and detailed level decisions. So the first two conditions are satisfied.  Then, for the last two conditions, the definitions equivalently imply that,
\begin{equation}
    \label{eq-proof-block-3}
    \begin{pmatrix}
            A^{h} & 0\\
            0 &  B^{h}
        \end{pmatrix} 
        \begin{pmatrix}
            {z}^{h}(\omega_{h})\\
            {y}^{h}(\omega_{h})
        \end{pmatrix} = \begin{pmatrix}
            b(\omega_{h})\\
            d(\omega_{h})
        \end{pmatrix} - \begin{pmatrix}
            R(\omega^{h-1}) & 0\\
            S(\omega^{h-1}) &  0
        \end{pmatrix}
        \begin{pmatrix}
            {z}^{h-1}\\
            {y}^{h-1}
        \end{pmatrix}.
\end{equation}
Then, we can derive,

\begin{equation}
\label{eq-proof-blockb}
    \left\{
\begin{array}{lr}
    A^{h}{z}^{h}(\omega_{h}) = b(\omega_{h})- R(\omega^{h-1}){z}^{h-1}& (a)\\
     B^{h}{y}^{h}(\omega_{h})= d(\omega_{h})-   S(\omega^{h-1})  {y}^{h-1}& (b)
\end{array}
\right.
\end{equation}

From Equation ~\ref{eq-proof-block-3}, two characteristics can be found for the problem with block-separable structure. First, the aggregate level decisions have a relationship between consecutive stages. Second, the detailed level decisions are only dependent on the aggregate level decisions from the same stage and have no direct effect on future decisions. In this paper, for each stage $h$, Constraints ~\ref{model3-stoch-flow conservation} - \ref{model3-stoch-num of nurses1} capture the relationship between aggregate level decisions from stage to stage, which are in alignment with Equation ~\ref{eq-proof-block-3}. Constraint ~\ref{model3-stoch-num of nurses2} defines that the staff availability decision (detailed level) has no direct effect on future constraints and is only dependent on the staffing decisions from the previous stage. The left constraints are all working within the same stage $h$, and are dependent on the detailed level decisions only. The nurse scheduling decisions for each stage show no direct effect on the staff sizing decisions in the following stage, which can be described by Equation ~\ref{eq-proof-block-3}. So \ref{model4-stoch-linear-obj} has a block-separable recourse.

\section{Proof of Remark 1}
\label{app-proof-remark-1}
% \proof{Proof of Remark 1.} 
\label{Proof-2stage}
First, if we define a dynamic programming type of recursion for the last stage $|\mathcal{H}|$, then we have,

\begin{align}
    \label{solution-milp-proof1}
   Q_{|\mathcal{H}|}\left(  {\boldsymbol{\chi}}_{|\mathcal{H}|-1}, {\boldsymbol{\xi}}_{|\mathcal{H}|}(\boldsymbol{\omega}_{|\mathcal{H}|})\right) =  minimize & \notag\\
   \quad
\overbrace{C_{U}\bar{n}_{|\mathcal{H}|}(\boldsymbol{\omega}_{|\mathcal{H}|}) + C_{O}\tilde{n}_{|\mathcal{H}|}(\boldsymbol{\omega}_{|\mathcal{H}|})}^{Staffing} \notag\\ 
+\overbrace{\sum\limits_{i\in\mathcal{I}}{\sum\limits_{d\in\mathcal{D}_{|\mathcal{H}|}}{\sum\limits_{s\in\mathcal{S}_{i}}{C_{A}{U}_{i,|\mathcal{H}|}^{d,s}(\boldsymbol{\omega}_{|\mathcal{H}|})}}}}^{Adjustments} \notag\\
+ \overbrace{\sum\limits_{d\in\mathcal{D}_{|\mathcal{H}|}}{\sum\limits_{j\in\mathcal{J}}{ C_{C}\left(\bar{q}_{j,|\mathcal{H}|}^{d}(\boldsymbol{\omega}_{|\mathcal{H}|})+\tilde{q}_{j,|\mathcal{H}|}^{d}(\boldsymbol{\omega}_{|\mathcal{H}|})\right) }}}^{Coverage} 
\notag\\ + 
            \overbrace{\sum\limits_{i\in\mathcal{I}}{
            \sum\limits_{m=1}^{V_{i,|\mathcal{H}|}}{C_{m}v_{i,|\mathcal{H}|}^{m}(\boldsymbol{\omega}_{|\mathcal{H}|})}
            }}^{Violations}  
\end{align}

Let $Q_{h+1}\left(  {\boldsymbol{\chi}}_{h}\right)=\mathbb{E}_{{\boldsymbol{\xi}}_{h+1}}\left[ Q_{h+1}\left(  {\boldsymbol{\chi}}_{h}, {\boldsymbol{\xi}}_{h+1}(\boldsymbol{\omega}_{h+1})\right)\right]$ for all $h = 2,..., |\mathcal{H}|-1$, the following recursion is obtained,

\begin{align}
    \label{solution-milp-proof2}
   Q_{h}\left(  {\boldsymbol{\chi}}_{h-1}, {\boldsymbol{\xi}}_{h}(\boldsymbol{\omega}_{h})\right) =  minimize & \notag\\
   \quad
\overbrace{C_{U}\bar{n}_{h}(\boldsymbol{\omega}_{h})+ C_{O}\tilde{n}_{h}(\boldsymbol{\omega}_{h}) }^{Staffing}
\notag\\ +\overbrace{\sum\limits_{i\in\mathcal{I}}{\sum\limits_{d\in\mathcal{D}_{h}}{\sum\limits_{s\in\mathcal{S}_{i}}{C_{A}U_{i,h}^{d,s}(\boldsymbol{\omega}_{h})}}}}^{Adjustments} 
\notag\\ + \overbrace{\sum\limits_{d\in\mathcal{D}_{h}}{\sum\limits_{j\in\mathcal{J}}{ C_{C}\left(\bar{q}_{j,h}^{d}(\boldsymbol{\omega}_{h})+\tilde{q}_{j,h}^{d}(\boldsymbol{\omega}_{h})\right) }}}^{Coverage} 
\notag\\ + 
            \overbrace{\sum\limits_{i\in\mathcal{I}}{
            \sum\limits_{m=1}^{V_{i,h}}{C_{m}v_{i,h}^{m}(\boldsymbol{\omega}_{h})}
            }}^{Violations} + Q_{h+1}\left(  {\boldsymbol{\chi}}_{h}\right)
\end{align}

Then the object of the problem is equivalent to the following,
\begin{equation}
    \label{solution-milp-proof3}
    minimize\quad 
    \overbrace{U_{0}\left( \boldsymbol{\chi} \right) }^{\text{stage 0 cost}}+ Q_{2}({\boldsymbol{\chi}}_{1})
\end{equation}
Due to the block-separable structure, we can separate the object into two parts,
\begin{multline}
    \label{solution-milp-proof4}
     Q_{h}\left(  {\boldsymbol{\chi}}_{h-1}, {\boldsymbol{\xi}}_{h}(\boldsymbol{\omega}_{h})\right) =  Q_{h}^{A}\left(  {\boldsymbol{\chi}}_{h-1}^{A}, {\boldsymbol{\xi}}_{h}(\boldsymbol{\omega}_{h})\right) \\ + 
     Q_{h}^{D}\left(  {\boldsymbol{\chi}}_{h-1}^{A}, {\boldsymbol{\xi}}_{h}(\boldsymbol{\omega}_{h})\right)
\end{multline}
where $ {\boldsymbol{\chi}}_{h}^{A}$ represents the set of all aggregate level decisions at stage $h$ and ${\boldsymbol{\chi}}_{h}^{D}$ is the set of all detailed level decisions at stage $h$. In this way, we aim to exclude the term ${\boldsymbol{\chi}}_{h-1}^{D}$ in ${\boldsymbol{\chi}}_{h-1}$. Specifically,

\begin{multline}
    \label{solution-milp-proof5}
  Q_{h}^{A}\left(  {\boldsymbol{\chi}}_{h-1}^{A}, {\boldsymbol{\xi}}_{h}(\boldsymbol{\omega}_{h})\right) =  minimize
  \quad
  Q_{h+1}\left(  {\boldsymbol{\chi}}_{h}\right)
  \\ + \overbrace{  C_{U}\bar{n}_{h}(\boldsymbol{\omega}_{h})+ C_{O}\tilde{n}_{h}(\boldsymbol{\omega}_{h}) }^{Staffing} 
\end{multline}

\begin{multline}
    \label{solution-milp-proof6}
   Q_{h}^{D}\left(  {\boldsymbol{\chi}}_{h-1}^{A}, {\boldsymbol{\xi}}_{h}(\boldsymbol{\omega}_{h})\right) =  minimize \\
 \overbrace{\sum\limits_{d\in\mathcal{D}_{h}}{\sum\limits_{j\in\mathcal{J}}{ C_{C}\left(\bar{q}_{j,h}^{d}(\boldsymbol{\omega}_{h})+\tilde{q}_{j,h}^{d}(\boldsymbol{\omega}_{h})\right) }}}^{Coverage} 
 \\ + 
            \overbrace{\sum\limits_{i\in\mathcal{I}}{
            \sum\limits_{m=1}^{V_{i,h}}{C_{m}v_{i,h}^{m}(\boldsymbol{\omega}_{h})}
            }}^{Violations}
\end{multline}

Then, the value we seek is shown in Equation \ref{solution-milp-proof7}.
\begin{figure*}
% \caption{Full objective}
\label{solution-milp-proof_fig}
\begin{align}
    \label{solution-milp-proof7}
 min \quad&\overbrace{U_{0}\left( \boldsymbol{\chi} \right) }^{\text{stage 0 cost}}+ Q_{2}({\boldsymbol{\chi}}_{1})\notag\\
&=\overbrace{U_{0}\left( \boldsymbol{\chi} \right) }^{\text{stage 0 cost}} +\mathbb{E}_{{\boldsymbol{\xi}}_{2}}\left[ Q_{2}\left(  {\boldsymbol{\chi}}_{1}, {\boldsymbol{\xi}}_{2}(\boldsymbol{\omega}_{2})\right)\right]\notag\\
& =\overbrace{U_{0}\left( \boldsymbol{\chi} \right) }^{\text{stage 0 cost}} +\mathbb{E}_{{\boldsymbol{\xi}}_{2}}\left[ Q_{h}^{A}\left(  {\boldsymbol{\chi}}_{1}^{A}, {\boldsymbol{\xi}}_{2}(\boldsymbol{\omega}_{2})\right) +   Q_{h}^{D}\left(  {\boldsymbol{\chi}}_{1}^{A}, {\boldsymbol{\xi}}_{2}(\boldsymbol{\omega}_{2})\right)\right]\notag\\
& =\overbrace{U_{0}\left( \boldsymbol{\chi} \right) }^{\text{stage 0 cost}} +\mathbb{E}_{{\boldsymbol{\xi}}_{2}} \left[ min\quad
\overbrace{ C_{U}\bar{n}_{2}(\boldsymbol{\omega}_{2})+ C_{O}\tilde{n}_{2}(\boldsymbol{\omega}_{2}) }^{Staffing} + Q_{3}\left(  {\boldsymbol{\chi}}_{2}\right)\right. \notag \\
&\quad\left.+min 
\overbrace{\sum\limits_{d\in\mathcal{D}_{2}}{\sum\limits_{j\in\mathcal{J}}{ C_{C}\left(\bar{q}_{j,2}^{d}(\boldsymbol{\omega}_{2})+\tilde{q}_{j,2}^{d}(\boldsymbol{\omega}_{2})\right) }}}^{Coverage} + 
            \overbrace{\sum\limits_{i\in\mathcal{I}}{
            \sum\limits_{m=1}^{V_{i,2}}{C_{m}v_{i,2}^{m}(\boldsymbol{\omega}_{2})}
            }}^{Violations}\right] \notag\\
& = \overbrace{U_{0}\left( \boldsymbol{\chi} \right) }^{\text{stage 0 cost}} +\mathbb{E}_{{\boldsymbol{\xi}}_{2}} \left[ min\quad
\overbrace{C_{U}\bar{n}_{2}(\boldsymbol{\omega}_{2})+ C_{O}\tilde{n}_{2}(\boldsymbol{\omega}_{2}) }^{Staffing} + ... \right.\notag\\
& \quad+ \left.\mathbb{E}_{{\boldsymbol{\xi}}_{|\mathcal{H}|}}[ min\quad
\overbrace{  C_{U}\bar{n}_{|\mathcal{H}|}(\boldsymbol{\omega}_{|\mathcal{H}|})+ C_{O}\tilde{n}_{|\mathcal{H}|}(\boldsymbol{\omega}_{|\mathcal{H}|}) }^{Staffing} ] ...\right]\notag \\
& \quad+ \mathbb{E}_{{\boldsymbol{\xi}}_{2}} \left[ min 
\overbrace{\sum\limits_{d\in\mathcal{D}_{2}}{\sum\limits_{j\in\mathcal{J}}{ C_{C}\left(\bar{q}_{j,2}^{d}(\boldsymbol{\omega}_{2})+\tilde{q}_{j,2}^{d}(\boldsymbol{\omega}_{2})\right) }}}^{Coverage} + 
            \overbrace{\sum\limits_{i\in\mathcal{I}}{
            \sum\limits_{m=1}^{V_{i,2}}{C_{m}v_{i,2}^{m}(\boldsymbol{\omega}_{2})}
            }}^{Violations}   + ... \right. \notag\\
& \quad\left.+ \mathbb{E}_{{\boldsymbol{\xi}}_{|\mathcal{H}|}}[ min 
\overbrace{\sum\limits_{d\in\mathcal{D}_{|\mathcal{H}|}}{\sum\limits_{j\in\mathcal{J}}{ C_{C}\left(\bar{q}_{j,|\mathcal{H}|}^{d}(\boldsymbol{\omega}_{|\mathcal{H}|})+\tilde{q}_{j,|\mathcal{H}|}^{d}(\boldsymbol{\omega}_{|\mathcal{H}|})\right) }}}^{Coverage} + 
            \overbrace{\sum\limits_{i\in\mathcal{I}}{
            \sum\limits_{m=1}^{V_{i,|\mathcal{H}|}}{C_{m}v_{i,|\mathcal{H}|}^{m}(\boldsymbol{\omega}_{|\mathcal{H}|})}
            }}^{Violations}   ] ...\right]\notag \\
\end{align}
\end{figure*}

According to Assumption \ref{Assumption-discrete scenarios} and Assumption \ref{Assumption-known scenario tree}, for each stage, the expectation term $\mathbb{E}_{\boldsymbol{\xi}_h}$ can be replaced by the summation of scenario probability multiply the corresponding objective terms of all scenarios. Then we have the formulation given in Objectives \ref{model5-2stage-obj} and \ref{model5-2stage-stoch-obj} (in the main text).

\newpage
\section{Proof of Proposition ~\ref{proposition-bybengio}}
\label{app-proof-prop-2}
% \proof{Proof of Proposition 2}
Note that in this appendix we just give the brief proof. Refer to \cite{Bengio2021a} for the detailed proof.
For any node in the flow network, the probability of reaching this node is as follows.
\begin{equation}
    \label{eq-appendixA-1}
    P\left(s^{'}\right)=\sum\limits_{(a,s):T(s,a)=s^{'}}{\pi\left(a|s\right)P(s)}.
\end{equation}
Together with Equation ~\ref{eq-proposition-begio1}, we have,
\begin{equation}
    \label{eq-appendixA-2}
    P\left(s^{'}\right)=\sum\limits_{(a,s):T(s,a)=s^{'}}{\frac{F(s,a)}{F(s)}P(s)}.
\end{equation}.
Next, the following Equation ~\ref{eq-appendixA-3} is proven by induction.
\begin{equation}
    \label{eq-appendixA-3}
    P\left(s\right)=\frac{F(s)}{F(s_{0})}.
\end{equation}
For the source node, this is trivially true since $P\left(s_{0}\right)=1$. For the following of nodes $s$ with successor $s^{'}$, if Equation ~\ref{eq-appendixA-3} is true, then
\begin{equation}
    \label{eq-appendixA-3b}
      \begin{aligned}
      P\left(s^{'}\right)&=\sum\limits_{(a,s):T(s,a)=s^{'}}{\frac{F(s,a)}{F(s)}\frac{F(s)}{F(s_{0})}}\\
      &=\frac{\sum\limits_{(a,s):T(s,a)=s^{'}}{F\left(s,a\right)}}{F(s_{0})}\\
      & = \frac{F(s^{'})}{F(s_{0})},
      \end{aligned}
\end{equation}
which is just  Equation ~\ref{eq-appendixA-3}.

Now we consider Equation ~\ref{eq-appendixA-3} with terminal states $x$, whose flow is $F(x)=R(x)$ and we have,
\begin{equation}
    \label{eq-appendixA-4}
    P(x)=\frac{R(x)}{F(s_{0})}.
\end{equation}

Note that $\sum_{x\in\mathcal{P}_{x}}{P(x)}=1$ and summing both sides in Equation~\ref{eq-appendixA-4}, there is,
\begin{equation}
    \label{eq-appendixA-4b}
    P(x)=\frac{R(x)}{\sum_{x^{'}\in\mathcal{P}_{x}}{R(x^{'})}}
\end{equation}

\newpage

\begin{turnpage}
\section{Summary of Comparison of Relevant Studies}
\label{appen-studies}
\begin{table}
\caption{\textbf{Relevant studies}}
\label{Table-literature review}
 {\def\arraystretch{1.2}
\begin{tabular}{m{1.5cm}|p{2.5cm}|p{2cm}|p{1.5cm}|p{2cm}|p{2cm}|p{2cm}|p{5cm}}
\hline\hline
\textbf{Studies} & \textbf{Demand settings} & \textbf{Uncertainty} & \textbf{Fairness} & \textbf{Flexibility} & \textbf{Regularity} & \textbf{Model} & \textbf{Solution method} \\ \hline
\cite{Ricardo2023} & deterministic  & $\checkmark$  & $\checkmark$   & $\checkmark$ & & MOP  & NSGAII and policy-based heuristics  \\
\cite{Anderson2023} &deterministic  &  &  & $\checkmark$ &  & HOM & hybrid exact and heuristic method\\
\cite{sun2023equitable} & stochastic &$\checkmark$ & $\checkmark$ & $\checkmark$ & & MIP & $\epsilon$-constraint solution method\\ 
\cite{zaerpour2022scheduling} & stochastic & $\checkmark$ & $\checkmark$ &  & & two-stage SP & L-shaped method\\
\cite{rath2022staff} & stochastic &$\checkmark$ & & $\checkmark$ & & two-stage SP & SAA-based heuristic\\
\cite{Guo2022} & deterministic &  &  & $\checkmark$ &  & MIP & column generation-based heuristic\\
\cite{Chen2022} &deterministic &  &  &  $\checkmark$ & $\checkmark$  & NN & heuristic strategy\\
\cite{Pieter2021} & deterministic  &  &  $\checkmark$ &  $\checkmark$ & $\checkmark$ &  MOP  & lexicographic goal programming with acceptance threshold  \\
\cite{Kheiri2021} & weekly demand   &  &  & $\checkmark$ &  & hidden MDP & local search and perturbation heuristics\\
\cite{Schoenfelder2020} & stochastic  & $\checkmark$  &  & $\checkmark$ &  & multi-stage SP &\\
\cite{Lai2020} & deterministic &   &  &  $\checkmark$ &  & NFM & graph construction algorithm and solver\\
\cite{Legrain2020} & weekly demand   &  $\checkmark$ &  & $\checkmark$ &  &  DP &online stochastic algorithm\\
\cite{Strandmark2020} & deterministic &   &  & $\checkmark$ &  &IP  &column generation-based heuristic\\
\cite{Wolbeck2020} & deterministic &   & $\checkmark$  &  &  & IP & commercial solver\\
\cite{Wickert2019} &  & nurse absence  &  &  &  & IP & variable neighborhood descent heuristic\\
\cite{Kim2015} & stochastic &  $\checkmark$ &  &  &  & two-stage SP & integer L-shaped algorithm\\
\cite{Wong2014} & deterministic &   &  & $\checkmark$ &  &two-stage model  & hybrid heuristic\\
\cite{Maenhout2013a} & deterministic  &   & $\checkmark$ & $\checkmark$ &  & IP & branch-and-price\\
\cite{Maenhout2013b} & deterministic &   & $\checkmark$ & $\checkmark$ &  &  IP&  branch-and-price and meta-heuristic\\
\cite{Wright2013} & deterministic &   &  & $\checkmark$  &  & bicriteria IP &commercial solver\\
\cite{Martin2013} & deterministic &  &  $\checkmark$ &  &  & agent-based model  &meta-heuristic\\
% This paper  & weekly demand   & weekly preference    &  $\checkmark$ &  $\checkmark$ &  $\checkmark$ &  $\checkmark$ &multi-stage SP & generative flow networks \\ 
\hline\hline
% \multicolumn{8}{l}{\scriptsize Note: MOP = multi-objective programming; HOM = hierarchical optimization model; NN = neural network; MDP = Markov decision programming; \\ NFM = network flow model; DP = dynamic programming; LP = linear programming; \\ MILP = mixed-integer programming; MIP = mixed-integer programming; SP = stochastic programming; IP = integer programming. }
\end{tabular}
}
\end{table}
\end{turnpage}

From Table \ref{Table-literature review}, we observe that most of the nurse scheduling problems studied in the literature still rely on the assumption that the problem is deterministic on patient demand. This may result in suboptimal solutions when fulfilling patient demand. Most existing studies considering uncertainties in nurse scheduling problems adopt a two-stage stochastic programming model. An implicit assumption is introduced that uncertain demand, though for multiple stages, is immediately revealed for the entire planning horizon. This assumption is more suitable for a short-term single-period system. However, in nurse scheduling problems, instead of knowing the uncertain information all at one time,  the realization of demand is more of a stage-by-stage manner. To enable this modeling capability, this paper develops a multi-stage stochastic programming model. 

%% file: 07_app_model1.tex
% \begin{figure*}
% \caption{Model 1 Section 1}
% \label{model1-section1}
\begin{widetext}
\section{Model 1: Deterministic Model}
\label{app-model1}
\begin{align}
    \label{model1-obj}
    minimize \quad U_{1}\left( \boldsymbol{\chi} \right) 
 &=\overbrace{C_{S}\sum\limits_{i\in\mathcal{I}}{\beta_{i}}}^{Staffing} + \overbrace{C_{C}\sum\limits_{j\in\mathcal{J}}{\sum\limits_{d\in\mathcal{D}}{(\bar{q}_{j}^{d}+\tilde{q}_{j}^{d})}}}^{Coverage} + \overbrace{
            C_{R}\sum\limits_{i\in\mathcal{I}}{\sum\limits_{d\in \mathcal{D}}{\sum\limits_{s\in\mathcal{S}_{i}}{R_{i}^{d,s}(1-x_{i}^{d,s})\beta_{i}}}}}^{Requests}  + 
            \overbrace{\sum\limits_{i\in\mathcal{I}}{
            \sum\limits_{m=1}^{V_{i}}{C_{m}v_{i}^{m}}
            }}^{Violations}
\end{align}

subject to,

\begin{equation}
    \label{model1-num of nurses2}
        \sum\limits_{i\in\mathcal{I}}{\beta_{i}} \leq N,  
\end{equation}

\begin{equation}
    \label{model1-h1}
        \sum\limits_{s\in \mathcal{S}_{i}}{x_{i}^{d,s}} \leq \beta_{i}, \quad \forall i \in \mathcal{I}, d \in \mathcal{D}
\end{equation}

\begin{equation}
    \label{model1-h2}
    T_{i}^{min}\beta_{i} \leq \sum\limits_{s\in \mathcal{S}_{i}}{\sum\limits_{d \in \mathcal{D}}{
    L_{s}x_{i}^{d,s}
    }}
    \leq T_{i}^{max}\beta_{i}, \quad \forall i \in \mathcal{I}
\end{equation}

% \textbf{H3. Work policy constraints:}
\begin{equation}
    \label{model1-h5}
    \sum\limits_{j\in \mathcal{J}}{y_{i}^{j}} \leq
    W_{p}, \quad \forall i \in \mathcal{I}_{p}, p\in\mathcal{P}
\end{equation}

\begin{equation}
    \label{model1-h6}
    y_{i}^{j} \geq 
    \sum\limits_{s \in \mathcal{S}_{i}\cap s \in\mathcal{S}_{j}}{x_{i}^{d,s}},\quad \forall i \in \mathcal{I},  d \in \mathcal{D}, j\in\mathcal{J}
\end{equation}

% \textbf{H4. ND-NW-AM avoided constraints:}
    \begin{align}
    \label{model1-h7}
    \left(1- \sum\limits_{s\in\mathcal{S}_{i}\cap \mathcal{S}_{N}}{
    x_{i}^{d,s}
    } \right) 
    +\sum\limits_{s\in\mathcal{S}_{i}}{
    x_{i}^{d+1,s}
    }
    +  \left(1- \sum\limits_{s\in\mathcal{S}_{i}\cap \mathcal{S}_{AM}}{
    x_{i}^{d+2,s}
    } \right)  \geq \beta_{i}, \quad
    \forall i \in \mathcal{I}, d \in \mathcal{D}\backslash D^{L}_{2}
\end{align}

% \textbf{H5. ND-AM/PM avoided constraints:}
\begin{align}
    \label{model1-h7b} \sum\limits_{s\in\mathcal{S}_{i}\cap \mathcal{S}_{N}}{
    x_{i}^{d,s}
    } + \sum\limits_{s\in\mathcal{S}_{i}:s\in \mathcal{S}_{AM}\cup \mathcal{S}_{PM}}{
    x_{i}^{d+1,s}
    } \leq \beta_{i}, \quad
    \forall i \in \mathcal{I}, d \in \mathcal{D}\backslash D^{L}_{1}
\end{align}

% \textbf{H6. Minimum nonworking constraints:}
\begin{equation}
    \label{model1-h8}
    7 - \sum\limits_{d\in\mathcal{D}_{k}}{
        \sum\limits_{s \in \mathcal{S}_{i}}{
        x_{i}^{d,s}
        }
    } \geq E_{i}\beta_{i}
    \quad \forall i \in \mathcal{I}, k\in \mathcal{K}
\end{equation}

% \textbf{H. Consecutive working days constraints:}
\begin{align}
    \label{model1-h16}
    \sum\limits_{d=d^{'}}^{d^{'}+F_{i}}{\sum\limits_{s\in\mathcal{S}_{i}}{x_{i}^{d,s}}}
    \leq F_{i}\beta_{i},\quad \forall i \in\mathcal{I}, d^{'} \in \mathcal{D}\backslash\mathcal{D}_{F_{i}-1}^{L}
\end{align}

% \textbf{H7. Demand coverage constraints:}
\begin{equation}
    \label{model1-h13}
    \sum\limits_{i\in\mathcal{I}}{
        \sum\limits_{s\in\mathcal{S}_{i}\cap \mathcal{S}_{j}}{
            x_{i}^{d,s}
        }
    } + \bar{q}_{j}^{d} - \tilde{q}_{j}^{d} 
    = Q_{j}^{d},
    \quad \forall d \in \mathcal{D}, j \in\mathcal{J}
\end{equation}

% \textbf{S5. Minimum nonworking weekend constraints:}
\begin{equation}
\label{model1-minNWweekends}
    \sum\limits_{d\in\mathcal{D}_{Sat}\cup \mathcal{D}_{Sun}}{\sum\limits_{s\in\mathcal{S}_{i}}{x_{i}^{d,s}}} 
  - p_{1,i}
   \leq M_{i}\beta_{i},
    \quad \forall i \in\mathcal{I}
\end{equation}

% \textbf{S1. Undesired pattern NW-AM constraints:}
\begin{align}
    \label{model1-h14}
   \sum\limits_{s\in\mathcal{S}_{i}}{x_{i}^{d,s}} + \left(1- 
  \sum\limits_{s\in\mathcal{S}_{i}\cap\mathcal{S}_{AM}}{x_{i}^{d+1,s}}   \right) +  p_{2,i}^{d} \geq \beta_{i},
    \quad \forall i \in\mathcal{I}, d \in \mathcal{D}\backslash D^{L}_{1}
\end{align}

% \end{figure*}

% \begin{figure*}
% \caption{Model 1 Section 2}
\label{model1-section2}

% \textbf{S6. Undesired pattern ND-ND on weekends constraints: }
\begin{align}
    \label{model1-h17}
    % Undesired to have two evening shifts on weekend (E-E)
    \left(1- \sum\limits_{s\in\mathcal{S}_{i}\cap \mathcal{S}_{N}}{
    x_{i}^{d,s}
    } \right) + \left(1-\sum\limits_{s\in\mathcal{S}_{i}\cap \mathcal{S}_{N}}{
    x_{i}^{d+1,s}
    } \right) + p_{3,i}^{d} \geq \beta_{i},\quad \forall i \in\mathcal{I}, d\in\mathcal{D}_{Sat}
\end{align}

% \textbf{S7. Undesired pattern NW-W-NW constraints:}
\begin{align}
    \label{model1-h18}
    % Undesired to have non-working before and after a working shift (NW-W-NW). Example: DO-D74-RD
    \sum\limits_{s\in\mathcal{S}_{i}}{x_{i}^{d,s}} + \left( 1- \sum\limits_{s\in\mathcal{S}_{i}}{x_{i}^{d+1,s}}\right) + \sum\limits_{s\in\mathcal{S}_{i}}{x_{i}^{d+2,s}} + p_{4,i}^{d} \geq \beta_{i},\quad\forall i \in\mathcal{I}, d\in\mathcal{D}\backslash D^{L}_{2}
\end{align}

\begin{equation}
    \label{model1-violations}
    \sum\limits_{m=1}^{V_{i}}{mv_{i}^{m}}= p_{1,i} + \sum\limits_{d \in \mathcal{D}\backslash D^{L}_{1}}{p_{2,i}^{d}} + \sum\limits_{d \in \mathcal{D}_{Sat}}{p_{3,i}^{d}} + \sum\limits_{d \in \mathcal{D}\backslash D^{L}_{2}}{p_{4,i}^{d}}, \quad \forall i \in\mathcal{I}
\end{equation}

% \textbf{Variables domain constraints:}
\begin{equation}
    \label{model1-variable domain_xy}
    \beta_{i}, x_{i}^{d,s}, y_{i}^{j} \in \{0,1\}, \quad \forall i \in\mathcal{I}, d\in\mathcal{D}, s\in \mathcal{S}_{i}, j\in\mathcal{J}
\end{equation}

\begin{equation}
    \label{model1-variable domain_g}
    \bar{q}_{j}^{d}, \tilde{q}_{j}^{d} \in \mathbb{Z}^{+}, \quad \forall j\in\mathcal{J}, d\in\mathcal{D}
\end{equation}

\begin{equation}
    \label{model1-variable domain_v}
    {v}_{i}^{m} \in \{0,1\}, \quad \forall i\in\mathcal{I}, m\in[0,V_{i}]
\end{equation}

\begin{equation}
    \label{model1-variable domain_p2}
    p_{1,i}  \in \mathbb{Z}^{+}\quad \forall i \in\mathcal{I}
\end{equation}

\begin{equation}
    \label{model1-variable domain_p1}
    p_{2,i}^d  \in \{0,1\}, \quad \forall i \in\mathcal{I}, d \in \mathcal{D}\backslash D^{L}_{1}
\end{equation}

\begin{equation}
    \label{model1-variable domain_p6}
    p_{3,i}^{d} \in \{0,1\}, \quad \forall i \in\mathcal{I}, d\in\mathcal{D}_{Sat}
\end{equation}

\begin{equation}
    \label{model1-variable domain_p3}
   p_{4,i}^{d} \in \{0,1\}, \quad \forall i \in\mathcal{I}, d \in \mathcal{D}\backslash D^{L}_{2}
\end{equation}

% \end{figure*}

\end{widetext}

%% file: 07_app_model2.tex
\begin{widetext}
\section{Model 2: Linearization of the Deterministic model}
\label{app-linear-deterministic}
\setlength{\jot}{0.3pt}

\begin{align}
    \label{model2-obj}
    minimize \quad U_{2}\left( \boldsymbol{\chi} \right) 
 =\overbrace{C_{S}\sum\limits_{i\in\mathcal{I}}{\beta_{i}}}^{Staffing} 
 + \overbrace{C_{C}\sum\limits_{j\in\mathcal{J}}{\sum\limits_{d\in\mathcal{D}}{(\bar{q}_{j}^{d}+\tilde{q}_{j}^{d})}}}^{Coverage}
 \notag\\+ \overbrace{
            C_{R}\sum\limits_{i\in\mathcal{I}}{\sum\limits_{d\in \mathcal{D}}{\sum\limits_{s\in\mathcal{S}_{i}}{R_{i}^{d,s}Y_{i}^{d,s}}}}}^{Requests} 
+ 
\overbrace{\sum\limits_{i\in\mathcal{I}}{
\sum\limits_{m=1}^{V_{i}}{C_{m}v_{i}^{m}}
}}^{Violations}
\end{align}

\begin{equation}
\label{model2-linearizations}
\begin{cases}
Y_{i}^{d,s} \leq 1 - x_{i}^{d,s} & \\ 
Y_{i}^{d,s} \leq \beta_{i} & \\ 
Y_{i}^{d,s} \geq \beta_{i} - x_{i}^{d,s} & \quad \forall i\in\mathcal{I}, s\in \mathcal{S}_{i}, d \in \mathcal{D}\\
Y_{i}^{d,s} \in \{0,1\} & 
\end{cases}
\end{equation}

\begin{equation}
    \label{model2-others}
    \text{Constraints} \quad ~\ref{model1-num of nurses2}-~\ref{model1-variable domain_p3}
\end{equation}

\end{widetext}

%% file: 07_app_model3.tex
\begin{widetext}
\section{Model 3: Stochastic Model}
\label{app-model3}
% \begin{figure*}
% \caption{Model 3 Section 1}
% \label{model3-section1}
\begin{align}
    \label{model3-stoch-obj}
    &U_{3}^{*} =  minimize\quad 
    \overbrace{U_{2}\left( \boldsymbol{\chi} \right) }^{\text{stage 0 cost}}+
    \overbrace{\mathbb{E}_{\boldsymbol{\xi}_{1}} \Bigg [\Upsilon_{1}\left( \boldsymbol{\chi},\boldsymbol{\xi}_{1}(\boldsymbol{\omega}_{1}) \right) + ... +  \mathbb{E}_{\boldsymbol{\xi}_{|\mathcal{H}|}}\left[\Upsilon_{|\mathcal{H}|}\left( \boldsymbol{\chi}_{|\mathcal{H}|-1},\boldsymbol{\xi}_{|\mathcal{H}|}(\boldsymbol{\omega}_{|\mathcal{H}|})\right)\right]...
    \Bigg]}^{\text{expected recourse of stage 1,..., stage }|\mathcal{H}|}
\end{align}

subject to

\begin{equation}
    \label{model3-stoch-stageh-obj}
    \text{Constraints} \quad ~\ref{model2-linearizations}-~\ref{model2-others}
\end{equation}
where for each stage $h$, we have the recourse function defined as follows,

\begin{align}
    \label{model3-stoch-obj2}
   \Upsilon_{h}\left(  {\boldsymbol{\chi}}_{h-1}, {\boldsymbol{\xi}}_{h}(\boldsymbol{\omega}_{h})\right) =  minimize &\quad
\overbrace{C_{U}\bar{n}_{h}(\boldsymbol{\omega}_{h})+ C_{O}\tilde{n}_{h}(\boldsymbol{\omega}_{h})}^{Staffing}+\overbrace{\sum\limits_{i\in\mathcal{I}}{\sum\limits_{d\in\mathcal{D}_{h}}{\sum\limits_{s\in\mathcal{S}_{i}}{C_{A}|x_{i}^{d,s}- {x}_{i,h}^{d,s}(\boldsymbol{\omega}_{h})|}}}}^{Adjustments}  \notag \\
&
+ \overbrace{\sum\limits_{d\in\mathcal{D}_{h}}{\sum\limits_{j\in\mathcal{J}}{ C_{C}\left(\bar{q}_{j,h}^{d}(\boldsymbol{\omega}_{h})+\tilde{q}_{j,h}^{d}(\boldsymbol{\omega}_{h})\right) }}}^{Coverage} + 
            \overbrace{\sum\limits_{i\in\mathcal{I}}{
            \sum\limits_{m=1}^{V_{i,h}}{C_{m}v_{i,h}^{m}(\boldsymbol{\omega}_{h})}
            }}^{Violations}  
\end{align}
subject to,

\begin{equation}
    \label{model3-stoch-flow conservation}
    n_{h}(\boldsymbol{\omega}_{h}) = n_{h-1}(\boldsymbol{\omega}_{h-1}) + \bar{n}_{h}(\boldsymbol{\omega}_{h}) - \tilde{n}_{h}(\boldsymbol{\omega}_{h}), \quad  \boldsymbol{\omega}_{h-1} =\mathbf{a}(\boldsymbol{\omega}_{h}),  h=2,...,|\mathcal{H}|
\end{equation}

\begin{equation}
    \label{model3-stoch-num of nurses1}
        \tilde{n}_{h}(\boldsymbol{\omega}_{h}) \le n_{h-1}(\boldsymbol{\omega}_{h-1}),  \quad \boldsymbol{\omega}_{h-1} =\mathbf{a}(\boldsymbol{\omega}_{h}),  h=2,...,|\mathcal{H}|
\end{equation}

\begin{equation}
    \label{model3-stoch-num of nurses2}
        \sum\limits_{i\in\mathcal{I}}{\beta_{i,h}(\boldsymbol{\omega}_{h})} \leq n_{h}(\boldsymbol{\omega}_{h}), \quad  h=1,...,|\mathcal{H}|  
\end{equation}

% \textbf{H1. Single assignment constraints:}
\begin{equation}
    \label{model3-stoch-h1}
    \sum\limits_{s\in \mathcal{S}_{i}}{ {x}_{i,h}^{d,s}(\boldsymbol{\omega}_{h})} \leq \beta_{i,h}(\boldsymbol{\omega}_{h}), \quad \forall i \in \mathcal{I}, d \in \mathcal{D}_{h}, h=1,...,|\mathcal{H}|
\end{equation}

\begin{equation}
    \label{model3-h2}
    T_{i,h}^{min}\beta_{i,h}(\boldsymbol{\omega}_{h}) \leq \sum\limits_{s\in \mathcal{S}_{i}}{\sum\limits_{d \in \mathcal{D}}{
    L_{s}x_{i,h}^{d,s}(\boldsymbol{\omega}_{h})
    }}
    \leq T_{i,h}^{max}\beta_{i,h}(\boldsymbol{\omega}_{h}), \quad \forall i \in \mathcal{I}, h=1,...,|\mathcal{H}|
\end{equation}

% \textbf{H3. Work policy constraints:}
\begin{equation}
    \label{model3-stoch-h6}
    {y}_{i}^{j} \geq 
    \sum\limits_{s \in \mathcal{S}_{i}\cap s \in\mathcal{S}_{j}}{ {x}_{i,h}^{d,s}(\boldsymbol{\omega}_{h})},\quad \forall i \in \mathcal{I},  d \in \mathcal{D}_{h}, j\in\mathcal{J}, h=1,...,|\mathcal{H}|
\end{equation}

% \textbf{H4. ND-NW-AM avoided constraints:}
\begin{align}
    \label{model3-stoch-h7}
    \left(1- \sum\limits_{s\in\mathcal{S}_{i}\cap \mathcal{S}_{N}}{
    x_{i,h}^{d,s}(\boldsymbol{\omega}_{h})
    } \right) 
    +\sum\limits_{s\in\mathcal{S}_{i}}{
    x_{i,h}^{d+1,s}(\boldsymbol{\omega}_{h})
    }
    +  \left(1- \sum\limits_{s\in\mathcal{S}_{i}\cap \mathcal{S}_{AM}}{
    x_{i,h}^{d+2,s}(\boldsymbol{\omega}_{h})
    } \right)  \geq \beta_{i,h}(\boldsymbol{\omega}_{h}) \quad \forall i \in \mathcal{I}, d \in \mathcal{D}_{h}\backslash D^{L}_{2}, h=1,...,|\mathcal{H}|
\end{align}

% \textbf{H5. ND-AM/PM avoided constraints:}
\begin{align}
    \label{model3-stoch-h8}
    \sum\limits_{s\in\mathcal{S}_{i}\cap \mathcal{S}_{N}}{
    x_{i,h}^{d,s}(\boldsymbol{\omega}_{h})
    }
    +  \sum\limits_{s\in\mathcal{S}_{i}:s\in \mathcal{S}_{AM}\cup\mathcal{S}_{PM}}{
    x_{i,h}^{d+1,s}(\boldsymbol{\omega}_{h})
    }   \leq \beta_{i,h}(\boldsymbol{\omega}_{h}), \quad \forall i \in \mathcal{I}, d \in \mathcal{D}_{h}\backslash D^{L}_{1}, h=1,...,|\mathcal{H}|
\end{align}

% \end{figure*}

% \begin{figure*}
% \label{model3-section2}
% \caption{Model 3 Section 2}

% \textbf{H6. Minimum nonworking constraints:}
\begin{equation}
    \label{model3-stoch-h9}
    7 - \sum\limits_{d\in\mathcal{D}_{h}}{
        \sum\limits_{s \in \mathcal{S}_{i}}{
         {x}_{i,h}^{d,s}(\boldsymbol{\omega}_{h})
        }
    } \geq E_{i,h}\beta_{i,h}(\boldsymbol{\omega}_{h}),
    \quad \forall i \in \mathcal{I}, h=1,...,|\mathcal{H}|
\end{equation}

%S4. Consecutive working days constraint:}
% \textbf{S4. Consecutive working days constraints:}
\begin{align}
    \label{model3-stoch-h16a}
    \sum\limits_{d=d^{'}}^{d^{'}+F_{i}}{\sum\limits_{s\in\mathcal{S}_{i}}{ {x}_{i,h}^{d,s}(\boldsymbol{\omega}_{h})}}
    \leq F_{i,h}\beta_{i,h}(\boldsymbol{\omega}_{h}),\quad\forall i \in\mathcal{I}, d^{'} \in \mathcal{D}_{h}\backslash\mathcal{D}_{F_{i,h}-1}^{h, L}, h=1,...,|\mathcal{H}|
\end{align}

% \textbf{H7. Demand coverage constraints:}
\begin{equation}
    \label{model3-stoch-h13a}
    \sum\limits_{i\in\mathcal{I}}{
        \sum\limits_{s\in\mathcal{S}_{i}\cap \mathcal{S}_{j} }{
             {x}_{i,h}^{d,s}(\boldsymbol{\omega}_{h}) +  \bar{q}_{j,h}^{d}(\boldsymbol{\omega}_{h})- \tilde{q}_{j,h}^{d}(\boldsymbol{\omega}_{h})
        }
    } 
    =  {Q}_{j,h}^{d}(\boldsymbol{\omega}_{h}), 
    \quad \forall d \in \mathcal{D}_{h}, j \in\mathcal{J}, h=1,...,|\mathcal{H}|
\end{equation}

% \textbf{S5. Minimum nonworking weekend constraints:}
\begin{equation}
\label{model3-stoch-minNWweekends}
    \sum\limits_{d\in\mathcal{D}_{h}: d\in\mathcal{D}_{Sat}\cup\mathcal{D}_{Sun} }{\sum\limits_{s\in\mathcal{S}_{i}}{ {x}_{i,h}^{d,s}}} 
  -  {p}_{1,i}^{h}(\boldsymbol{\omega}_{h})
   \leq  M_{i,h}\beta_{i,h}(\boldsymbol{\omega}_{h}),
    \quad \forall i \in\mathcal{I}, h=1,...,|\mathcal{H}|
\end{equation}

% S1. Undesired pattern NW-AM constraint:
% \textbf{S1. Undesired pattern NW-AM constraints:}
\begin{align}
    \label{model3-stoch-h14}
   \sum\limits_{s\in\mathcal{S}_{i}}{ {x}_{i,h}^{d,s}(\boldsymbol{\omega}_{h})} + \left(1- 
  \sum\limits_{s\in\mathcal{S}_{i}\cap\mathcal{S}_{AM}}{ {x}_{i,h}^{d+1,s}(\boldsymbol{\omega}_{h})}   \right) +   {p}_{2,i}^{d,h}(\boldsymbol{\omega}_{h}) \geq \beta_{i,h}(\boldsymbol{\omega}_{h}),\quad \forall i \in\mathcal{I}, d \in \mathcal{D}_{h}\backslash\mathcal{D}^{h,L}_{1}, h=1,...,|\mathcal{H}| 
\end{align}

% \textbf{S6. Undesired pattern ND-ND on weekends: }
\begin{align}
    \label{model3-stoch-sh17}
    % Undesired to have two evening shifts on weekend (E-E)
    \left(1- \sum\limits_{s\in\mathcal{S}_{i}\cap \mathcal{S}_{N}}{
     {x}_{i,h}^{d,s}
    } \right) + \left(1-\sum\limits_{s\in\mathcal{S}_{i}\cap \mathcal{S}_{N}}{
     {x}_{i,h}^{d+1,s}
    } \right) +  {p}_{3,i}^{d,h}(\boldsymbol{\omega}_{h}) \geq \beta_{i,h}(\boldsymbol{\omega}_{h}),\quad\forall i \in\mathcal{I}, d\in\mathcal{D}_{Sat}\cap\mathcal{D}_{h}, h=1,...,|\mathcal{H}|
\end{align}

% \textbf{S7. Undesired pattern NW-W-NW constraints:}
%S7. Undesired pattern NW-W-NW constraint:}
\begin{align}
    \label{model3-stoch-h18}
    % Undesired to have non-working before and after a working shift (NW-W-NW). Example: DO-D74-RD
    \sum\limits_{s\in\mathcal{S}_{i}}{ {x}_{i,h}^{d,s}(\boldsymbol{\omega}_{h})} + \left( 1- \sum\limits_{s\in\mathcal{S}_{i}}{ {x}_{i,h}^{d+1,s}(\boldsymbol{\omega}_{h})}\right) + \sum\limits_{s\in\mathcal{S}_{i}}{ {x}_{i,h}^{d+2,s}(\boldsymbol{\omega}_{h})} +  {p}_{4,i}^{d,h}(\boldsymbol{\omega}_{h}) \geq \beta_{i}(\boldsymbol{\omega}_{h}),\quad\forall i \in\mathcal{I}, d\in\mathcal{D}_{h}\backslash\mathcal{D}^{h,L}_{2}, h=1,...,|\mathcal{H}|
\end{align}

\begin{equation}
    \label{model3-violations}
    \sum\limits_{m=1}^{V_{i,h}}{mv_{i,h}^{m}(\boldsymbol{\omega}_{h})}=  p_{1,i}^{h}(\boldsymbol{\omega}_{h})+\sum\limits_{d \in \mathcal{D}_{h}\backslash D^{L}_{1}}{p_{2,i}^{d,h}(\boldsymbol{\omega}_{h})}  + \sum\limits_{d \in \mathcal{D}_{h}\cap\mathcal{D}_{Sat}}{p_{3,i}^{d,h}(\boldsymbol{\omega}_{h})} + \sum\limits_{d \in \mathcal{D}_{h}\backslash D^{L}_{2}}{p_{4,i}^{d}(\boldsymbol{\omega}_{h})}, \quad \forall i \in\mathcal{I}, h=1,...,|\mathcal{H}|
\end{equation}
% \end{figure*}

% \begin{figure*}
% \label{model3-section3}
% \caption{Model 3 Section 3}
% \textbf{Variables domain constraints:}
\begin{equation}
    \label{model3-stoch-variable domain n}
    {n}_{h}(\boldsymbol{\omega}_{h}), \bar{n}_{h}(\boldsymbol{\omega}_{h}), \tilde{n}_{h}(\boldsymbol{\omega}_{h}) \in \mathbb{Z}^{+},   \quad  h=1,...,|\mathcal{H}|
\end{equation}

\begin{equation}
    \label{model3-stoch-variable domain_xy}
    \beta_{i,h}(\boldsymbol{\omega}_{h}), x_{i,h}^{d,s}(\boldsymbol{\omega}_{h}) \in \{0,1\}, \quad \forall i \in\mathcal{I}, d\in\mathcal{D}, s\in \mathcal{S}_{i}, j\in\mathcal{J},  h=1,...,|\mathcal{H}|
\end{equation}

\begin{equation}
    \label{model3-stoch-variable domain_g}
    \bar{q}^{d}_{j,h}(\boldsymbol{\omega}_{h}), \tilde{q}^{d}_{j,h}(\boldsymbol{\omega}_{h})\in \mathbb{Z}^{+}, \quad \forall j\in\mathcal{J}, d\in\mathcal{D}_{h},  h=1,...,|\mathcal{H}|
\end{equation}

\begin{equation}
    \label{model3-variable domain_p2}
    p_{1,i}^{h}(\boldsymbol{\omega}_{h})  \in \mathbb{Z}^{+}, \quad \forall i \in\mathcal{I},  h=1,...,|\mathcal{H}|
\end{equation}

\begin{equation}
    \label{model3-variable domain_p1}
    p_{2,i}^{d,h}(\boldsymbol{\omega}_{h})  \in \{0,1\}, \quad \forall i \in\mathcal{I}, d \in \mathcal{D}_{h}\backslash D^{L}_{1},  h=1,...,|\mathcal{H}|
\end{equation}

\begin{equation}
    \label{model3-variable domain_p6}
    p_{3,i}^{d,h}(\boldsymbol{\omega}_{h}) \in \{0,1\}, \quad \forall i \in\mathcal{I}, d\in\mathcal{D}_{Sat}\cap\mathcal{D}_{h},  h=1,...,|\mathcal{H}|
\end{equation}

\begin{equation}
    \label{model3-variable domain_p3}
   p_{4,i}^{d,h}(\boldsymbol{\omega}_{h})\in \{0,1\}, \quad \forall i \in\mathcal{I}, d \in \mathcal{D}_{h}\backslash D^{L}_{2},  h=1,...,|\mathcal{H}|
\end{equation}

% \end{figure*}
\end{widetext}

%% file: 07_app_model4.tex
\begin{widetext}
\section{Model 4: Linearization of the Stochastic model}
\label{app-linear-stoch}
Model 3 can be linearized by replacing the absolute term with one auxiliary variable $U_{i,h}^{d,s}(\boldsymbol{\omega}_{h})$ and the following additional constraints:
\begin{equation}
    \label{model4-stoch-linearize abs}
   x_{i}^{d,s}- {x}_{i,h}^{d,s}(\boldsymbol{\omega}_{h}) \leq U_{i,h}^{d,s}(\boldsymbol{\omega}_{h})
\end{equation}
\begin{equation}
    \label{model4-stoch-linearize abs2}
    {x}_{i,h}^{d,s}(\boldsymbol{\omega}_{h}) -x_{i}^{d,s} \leq U_{i,h}^{d,s}(\boldsymbol{\omega}_{h})
\end{equation}
Then the multi-stage stochastic nonlinear programming model can be equivalently linearized as shown in Equations \ref{model4-stoch-linear-obj} to \ref{model4-stoch-linear-stageh-constraints-3}.

% [Model 4]
% \begin{figure*}
% \label{Model-stoch-linearized}
% \caption{Model 4}
\begin{align}
    \label{model4-stoch-linear-obj}
    U_{4}^{*} &=minimize\quad 
    \overbrace{U_{2}\left( \boldsymbol{\chi} \right) }^{\text{stage 0 cost}}+
    \overbrace{\mathbb{E}_{\boldsymbol{\xi}_{1}} \Bigg [\Upsilon_{1}\left( \boldsymbol{\chi},\boldsymbol{\xi}_{1}(\boldsymbol{\omega}_{1}) \right) + ... +  \mathbb{E}_{\boldsymbol{\xi}_{|\mathcal{H}|}}\left[\Upsilon_{|\mathcal{H}|}\left( \boldsymbol{\chi}_{|\mathcal{H}|-1},\boldsymbol{\xi}_{|\mathcal{H}|}(\boldsymbol{\omega}_{|\mathcal{H}|})\right)\right]...
    \Bigg]}^{\text{expected recourse of stage 1,..., stage }|\mathcal{H}|} \\\notag
\end{align}

subject to
\begin{equation}
    \label{model4-stoch-stageh-constraints}
    \text{Constraints} \quad ~\ref{model2-linearizations}-~\ref{model2-others}
\end{equation}
where for each stage $h$, we have the recourse function defined as follows,
\begin{align}
    \label{model4-stoch-obj}
   \Upsilon_{h}\left(  {\boldsymbol{\chi}}_{h-1}, {\boldsymbol{\xi}}_{h}(\boldsymbol{\omega}_{h})\right) =  minimize &\quad
\overbrace{C_{U}\bar{n}_{h}(\boldsymbol{\omega}_{h})+ C_{O}\tilde{n}_{h}(\boldsymbol{\omega}_{h})}^{Staffing}+\overbrace{\sum\limits_{i\in\mathcal{I}}{\sum\limits_{d\in\mathcal{D}_{h}}{\sum\limits_{s\in\mathcal{S}_{i}}{C_{A}U_{i,h}^{d,s}(\boldsymbol{\omega}_{h})}}}}^{Adjustments} \notag \\
&
+ \overbrace{\sum\limits_{d\in\mathcal{D}_{h}}{\sum\limits_{j\in\mathcal{J}}{ C_{C}\left(\bar{q}_{h}^{d,j}(\boldsymbol{\omega}_{h})+\tilde{q}_{h}^{d,j}(\boldsymbol{\omega}_{h})\right) }}}^{Coverage} + 
            \overbrace{\sum\limits_{i\in\mathcal{I}}{
            \sum\limits_{m=1}^{V_{i,h}}{C_{m}v_{i,h}^{m}(\boldsymbol{\omega}_{h})}
            }}^{Violations} 
\end{align}

subject to,
\begin{equation}
    \label{model4-stoch-linear-stageh-constraints-2}
    \text{Constraints} \quad~\ref{model3-stoch-flow conservation}-~\ref{model3-variable domain_p3}
\end{equation}

\begin{equation}
    \label{model4-stoch-linear-stageh-constraints-3}
    \text{Constraints} \quad~\ref{model4-stoch-linearize abs}-~\ref{model4-stoch-linearize abs2}
\end{equation}
% \end{figure*}
\end{widetext}

%% file: 07_app_model5.tex
-
\clearpage
\begin{widetext}
\section{Model 5: Two-stage Stochastic Program with Recourse model}
\label{app-model5}
% \caption{Model 5}
\begin{align}
    \label{model5-2stage-obj}
    U_{5}^{*} = minimize \quad 
    &U_{2}\left( \boldsymbol{\chi} \right)  + \sum\limits_{h\in \mathcal{H}}{\sum\limits_{\boldsymbol{\omega}_{h} \in \mathcal{W}_{h} }{\boldsymbol{p}(\boldsymbol{\omega}_{h} ) \left [ C_{U}\bar{n}_{h}(\boldsymbol{\omega}_{h})+ C_{O}\tilde{n}_{h}(\boldsymbol{\omega}_{h}) \right]  }}\notag \\
    & + \sum\limits_{h\in \mathcal{H}}{\sum\limits_{\boldsymbol{\omega}_{h} \in \mathcal{W}_{h}}{  \sum\limits_{\boldsymbol{\omega}_{h+1}\in\mathcal{G}(\boldsymbol{\omega}_{h})}{\boldsymbol{p}(\boldsymbol{\omega}_{h+1})f_{h+1}\left(\boldsymbol{\chi},n_{h}(\boldsymbol{\omega}_{h}), \bar{n}_{h}(\boldsymbol{\omega}_{h}), \tilde{n}_{h}(\boldsymbol{\omega}_{h})\right) }}}
\end{align}
where

\begin{align}
    \label{model5-2stage-stoch-obj}
   f_{h+1}\left(\boldsymbol{\chi},n_{h}(\boldsymbol{\omega}_{h}), \bar{n}_{h}(\boldsymbol{\omega}_{h}), \tilde{n}_{h}(\boldsymbol{\omega}_{h})\right)& =  minimize\quad\overbrace{\sum\limits_{i\in\mathcal{I}}{\sum\limits_{d\in\mathcal{D}_{h+1}}{\sum\limits_{s\in\mathcal{S}_{i}}{C_{A}U_{i,h+1}^{d,s}(\boldsymbol{\omega}_{h+1})}}}}^{Adjustments}  \notag\\ & + \overbrace{\sum\limits_{d\in\mathcal{D}_{h+1}}{\sum\limits_{j\in\mathcal{J}}{ C_{C}\left(\bar{q}_{j,h+1}^{d}(\boldsymbol{\omega}_{h+1})+\tilde{q}_{j,h+1}^{d}(\boldsymbol{\omega}_{h+1})\right) }}}^{Coverage} \notag\\
   & + 
    \overbrace{\sum\limits_{i\in\mathcal{I}}{
    \sum\limits_{m=1}^{V_{i,h}}{C_{m}v_{i,h+1}^{m}(\boldsymbol{\omega}_{h+1})}
    }}^{Violations}   
\end{align}

subject to,
\begin{equation}
    \label{model5-solution-equi-milp-constraint}
    \text{Constraints} \quad~\ref{model4-stoch-stageh-constraints}, \ref{model4-stoch-linear-stageh-constraints-2}, \text{and } \ref{model4-stoch-linear-stageh-constraints-3}
\end{equation}

\end{widetext}